\documentclass{article}
\usepackage{newclude}
\usepackage{comment}
\usepackage{icomma}
\usepackage{amsmath}
\usepackage{color}
\usepackage{cases}
\usepackage{tikz-cd}
\usepackage{dsfont}
\usepackage{geometry}
\usepackage{amsthm}
\usepackage{aliascnt}
\usepackage{amssymb}
\usepackage{graphicx}
\usepackage{adjustbox}
\usepackage{mathrsfs}
\usepackage[
    backend=biber,
    style=numeric,
    sortlocale=de_DE,
    natbib=true,
    url=false, 
    doi=true,
    eprint=true
]{biblatex}
\usepackage{hyperref}
\usepackage{cleveref}
\numberwithin{equation}{section}

\theoremstyle{plain}
\newtheorem{theorem}{Theorem}[section]
\newaliascnt{lemma}{theorem}
\newtheorem{lemma}[lemma]{Lemma}
\aliascntresetthe{lemma}
\newaliascnt{proposition}{theorem}
\newtheorem{proposition}[proposition]{Proposition}
\aliascntresetthe{proposition}
\newaliascnt{corollary}{theorem}
\newtheorem{corollary}[corollary]{Corollary}
\aliascntresetthe{corollary}
\newaliascnt{conjecture}{theorem}

\aliascntresetthe{conjecture}

\theoremstyle{definition}
\newtheorem{definition}{Definition}[section]
\newtheorem{remark}{Remark}[section]
\newtheorem{example}{Example}[section]

\newtheorem{convention}{Convention}[section]

\theoremstyle{remark}   

\def\Z{\mathbb{Z}}

\def\C{\mathbb{C}}
\def\R{\mathbb{R}}

\def\T{\mathbb{T}}

\DeclareMathOperator{\id}{Id}
\DeclareMathOperator{\im}{im}

\DeclareMathOperator{\colim}{colim}
\DeclareMathOperator{\hocolim}{hocolim}

\def\scrl{\mathscr{L}}

\author{Yin Wei \\
School of Mathematics, Nanjing University, Nanjing 210093, P.R. China \\
\texttt{weiyin919@qq.com}}
\title{Equivariant cohomology rings and their $H^*(B\T^n;\Z)$-module structures of locally standard $\T^n$-manifolds}

\begin{document}
    \maketitle
    \begin{abstract}
        We compute the integral equivariant cohomology ring of any locally standard $\T^n$-manifold, together with its module structure over $H^*(B\T^n;\Z)$. The ring is isomorphic to the topological face ring of the orbit space. The module structure is determined by the characteristic classes of the principal $\T^n$-bundle determined by the free orbits.
    \end{abstract}
    \section{Introduction}
A $\T^n$-manifold means a smooth manifold equipped with a smooth $\T^n$-action. Let $\T^n$ denote the $n$-dimensional torus group. A $\T^n$-manifold $M$ is called \emph{locally standard} if every point admits an invariant neighborhood equivariantly diffeomorphic, up to an automorphism of $\T^n$, to an invariant open subset of $\C^k\times\T^{n-k}\times\R^r$, for some $0\leq k\leq n$ and $r\geq0$, see \cite[Section~1]{wiemeler2022smooth}. Here $\T^n$ is identified with $\T^k\times\T^{n-k}$, where the first factor acts on $\C^k$ by coordinatewise complex multiplication, the second acts on $\T^{n-k}$ by multiplication, and the action on $\R^r$ is trivial. 
The orbit space $M/\T^n$ is then a \emph{manifold with corners} (see \Cref{def:nicemfdwithcorner}).

Classical examples include quasitoric manifolds, whose orbit spaces are simple convex polytopes. These were studied extensively by Davis and Januszkiewicz in \cite{MR1104531} (where they are called \emph{toric manifolds}). In particular, they computed both the equivariant cohomology ring of such a manifold and its module structure over $H^*(B\T^n;\Z)$, and they showed that this ring is the \emph{face ring} of the convex polytope. As the name
indicates, its ring structure is completely determined by the combinatorics of the faces of the orbit space.

In \cite{MR2283418}, Masuda and Panov studied \emph{torus manifolds} and related their cohomological properties to the combinatorics of the orbit space. In particular, they showed that the odd-degree cohomology of a torus manifold vanishes if and only if the action is locally standard and the orbit space is \emph{face-acyclic}, in which case the equivariant cohomology ring is the face ring of the simplicial poset of faces of the orbit space. Yu solved the case where the principal bundle determined by the free orbit is trivial \cite{MR4721373}, showing that the equivariant cohomology ring is the \emph{topological face ring} (see \Cref{def:topologicalfacering}), a generalization of the face ring. Like the face ring, it is determined by data on the orbit space. Beyond the combinatorics recorded by the face ring, it also involves the cohomology rings of the faces and the restriction maps between them.

Xiong and Zeng \cite{xiong2023equivariantcohomologyringopen} computed the equivariant cohomology ring and its $H^*(B\T^n;\Z)$-module structure for locally standard $\T^n$-actions on open $2n$-dimensional torus manifolds, allowing nontrivial principal orbit bundles. 
Their proof proceeds by induction on the number of faces of the orbit space, with repeated use of the Mayer--Vietoris sequence. 
We treat the more general smooth locally standard actions defined above, allowing the manifold dimension to differ from $2n$ (see \Cref{THM:2}). 
We also prove an inverse-limit formula for the equivariant cohomology of the generalized \emph{polyhedral product}, which carries a $\T^{n+m}$-action (\Cref{thm:generalpolyhedral}).
Our proof uses a different method based on homotopy colimits and the vanishing of higher limits to compute both the ring and its module structure.

Our first main result computes the equivariant cohomology of the fiber product of a \emph{generalized moment-angle manifold} with a principal torus bundle (\Cref{THM:1}). Now we collect the definitions needed to make the statement precise.

We say that a manifold with corners is \emph{nice} if each of its codimension-$k$ faces is a connected component of the intersection of $k$ distinct facets. Throughout, nice manifolds with corners are assumed to be smooth and to have finitely many facets.
Let $Q$ be such a manifold with facet set $\mathcal F=\{F_1, \ldots, F_m\}$, and put $[m]=\{1, \ldots, m\}$. For $J\subseteq[m]$, set
\begin{equation}\label{eq:defoffacebycap}
    F_{\cap J}:=\bigcap_{j\in J}F_j, \qquad F_{\cap\varnothing}:=Q.
\end{equation}
For $J'\subseteq J$ we have $F_{\cap J}\subseteq F_{\cap J'}$. For a coefficient ring $R$ and $\alpha\in H^*(F_{\cap J'};R)$, write $\alpha|_{F_{\cap J}}$ for its pullback along this inclusion. 

We next introduce the polynomial notation used in the topological face ring.
Let $R$ be a commutative ring. For $J=\{j_1, \ldots, j_s\}$, we write
\[
    R[\mathbf x_J]:=R[x_{j_1}, \ldots, x_{j_s}]
\]
for the polynomial ring in the indicated variables. We use $R[\mathbf x_J^+]$ to denote its $R$-submodule spanned by monomials with strictly positive exponent in every variable:
\[
    R[\mathbf x_J^+]:=\operatorname{span}_R\bigl\{x_{j_1}^{n_1}\cdots x_{j_s}^{n_s}\mid n_1, \ldots, n_s\in\Z_{>0}\bigr\}.
\]
For $J\ne\varnothing$, this is a subring without unit. We set $R[\mathbf x_\varnothing]=R[\mathbf x_\varnothing^+]=R$. The topological face ring was introduced by Yu \cite{MR4721373}. We recall its definition here.

\begin{definition}\label{def:topologicalfacering}
    The \emph{topological face ring} of $Q$ over a graded commutative coefficient ring $R$ is
    \begin{equation}
        R\langle Q\rangle:=\bigoplus_{J\subseteq[m]}H^*(F_{\cap J};R)\otimes_R R[\mathbf x_J^+],
    \end{equation}
    where we take the convention that $H^*(\varnothing;R)=0$. Each $x_j$ has degree $2$, and the elements in $H^*(F_{\cap J};R)$, which are called \emph{coefficient classes}, have their canonical degrees as cohomology elements. For elements supported on $J$ and $J'$, multiplication is defined by
    \[
        (\alpha\otimes f)\star(\alpha'\otimes f')
        :=\bigl(\alpha|_{F_{\cap(J\cup J')}}\cup\alpha'|_{F_{\cap(J\cup J')}}\bigr)\otimes(ff').
    \]
\end{definition}
Thus multiplication restricts the coefficient classes to the common face intersection and takes their cup product. 

\begin{convention}\label{conv:main}
We collect three conventions that are used throughout.
\begin{enumerate}
    \item A parenthesized subscript denotes an indexed copy of a space, bundle or group. For instance, we write
    \[
        \T^m\cong\prod_{i\in[m]}\T^1_{(i)},
    \]
    so that $\T^1_{(i)}$ is the $i$-th circle factor of $\T^m$, which is a copy of $\T^1$. For $J\subseteq[m]$, we write
    \[
        \T^m_J:=\prod_{j\in J}\T^1_{(j)}\times\prod_{j\notin J}1_{(j)},
    \]
    the subgroup of $\T^m$ whose $j$-th coordinate is trivial for $j\notin J$.
    \item For a $G$-space $X$, we write $X_G:=EG\times_G X$ for its Borel construction. We use left actions on principal bundles and on $EG$, and take the diagonal quotient $(e, x)\sim(ge, gx)$ in the Borel construction.
    \item For a principal $\T^m$-bundle $\xi\to Q$ and a left $\T^m$-space $F$, we write
    \[
        \xi\times_{\T^m}F:=(\xi\times F)/{\sim},
        \qquad (v, z)\sim(gv, gz)\quad (g\in\T^m),
    \]
    for the associated bundle of $\xi$ with fiber $F$; thus $(v, z)$ is identified with its translate by the diagonal action of $\T^m$ on $\xi\times F$.
\end{enumerate}
\end{convention}

A \emph{$k$-dimensional unimodular subspace} of $\Z^l$ is a rank-$k$ submodule which is a direct summand. A function $\lambda:\mathcal F\to\Z^l$ is called a \emph{characteristic function} if, whenever $\varnothing\ne J\subseteq[m]$ and $F_{\cap J}\ne\varnothing$, the vectors $\{\lambda(F_j)\}_{j\in J}$ span a $|J|$-dimensional unimodular subspace.
Since $\lambda$ sends the $j$-th standard basis vector of $\Z^m$ to $\lambda(F_j)$, it induces a homomorphism, again denoted by $\lambda$, from $\Z^m$ to $\Z^l$, and hence a homomorphism of tori
\[
    \Lambda:\T^m\longrightarrow\T^l.
\]
Writing $g=(e^{2\pi i t_1},\ldots,e^{2\pi i t_m})\in\T^m$ and $\lambda(F_j)=(\lambda_{1j},\ldots,\lambda_{lj})$, we have
\[
    \Lambda(g)=\Bigl(e^{2\pi i\sum_{j=1}^m\lambda_{1j}t_j},\ \ldots,\ e^{2\pi i\sum_{j=1}^m\lambda_{lj}t_j}\Bigr).
\]
For $J\subseteq[m]$, we set $\T^l_{\lambda,J}:=\Lambda(\T^m_J)$, where $\T^m_J$ is the subgroup of \Cref{conv:main}.

\medskip
Let $E\to B$ be a principal $\T^r$-bundle, with a fixed identification of its structure group with $\T^r$. Write
\[
    c(E):=\bigl(c(E)_1, \ldots, c(E)_r\bigr)\in H^2(B;\Z^r)
\]
for its characteristic class, where the components are taken with respect to the coordinates of $\T^r$. Each component $c(E)_i$ lies in $H^2(B;\Z)$. For a principal circle bundle we identify $c(E)$ with its single component.

We now describe the general quotient construction of \emph{generalized moment-angle manifolds}, which is also used for the locally standard $\T^n$-manifolds of \Cref{sec:torus-manifolds}.
\begin{definition}\label{def:quotientconstruction}
    Let $\lambda:\mathcal F\to\Z^l$ be a characteristic function.
    For a principal $\T^l$-bundle $\pi_\xi:\xi\to Q$, define
\begin{equation}\label{eq:quotientconstruction}
    M_{(Q, \xi, \lambda)}:=\xi/\sim,
\end{equation}
where, for the unique $J\subseteq[m]$ with $\pi_\xi(p_1)\in F^\circ_{\cap J}$,
\begin{equation}\label{eq:relationdefquotient}
    p_1\sim p_2\ \Longleftrightarrow
    \pi_\xi(p_1)=\pi_\xi(p_2)\in F_{\cap J}^{\circ}
    \text{ and }p_2=gp_1\text{ for some }g\in\T^l_{\lambda, J}.
\end{equation}
Here $F_{\cap J}^{\circ}:=F_{\cap J}\setminus\bigcup_{i\notin J}F_i$, which is the relative interior of $F_{\cap J}$. 
\end{definition}

\begin{definition}\label{def:fiberproduct}
    Let $\pi_i:E_i\to Q$ $(i=1, 2)$ be maps from spaces $E_i$ to a common space $Q$. Their \emph{fiber product} over $Q$ is
    \begin{equation}\label{eq:fiberproduct}
        E_1\times_Q E_2:=\bigl\{(e_1, e_2)\in E_1\times E_2\mid \pi_1(e_1)=\pi_2(e_2)\bigr\}.
    \end{equation}
\end{definition}

Consider the fiber product of a principal $\T^n$-bundle $\zeta \to Q$ with a trivial principal $\T^m$-bundle $Q\times \T^m \to Q$. Their fiber product is a principal $\T^{n+m}$-bundle over $Q$:
\begin{equation*}
    \zeta \times_Q (Q \times \T^m) = \zeta \times \T^m
\end{equation*}
We define the generalized moment-angle manifold as a special case of the quotient construction given in \Cref{def:quotientconstruction}.
\begin{definition}
    In \Cref{def:quotientconstruction}, take $\xi = \zeta \times \T^m$ and $l=n+m$. 
    Let $\lambda(F_i):=e_{n+i}$ be the $n+i$-th standard basis vector of $\Z^{n+m}$, then the quotient space 
    \begin{equation}
        (\zeta \times \T^m) / \sim
    \end{equation}
    is called the \emph{generalized moment-angle manifold}, and denoted by $\mathscr L_{(Q, \zeta)}$.
\end{definition}

If $\zeta$ is removed from this construction, that is, if one takes $\xi=Q\times\T^m$ and $l=m$, the same quotient construction gives the generalized moment-angle manifold for the trivial $\T^m$-bundle over $Q$. Following Yu \cite{MR4721373}, we denote it by $\mathscr L_Q$.

By construction, $\mathscr L_{(Q,\zeta)}$ can also be written as the fiber product
\begin{equation*}
    \mathscr L_{(Q,\zeta)}\cong\zeta\times_Q\mathscr L_Q
\end{equation*}
over $Q$.

\begin{remark}\label{rem:fiberproductoverQ}
    Unless otherwise stated, the fiber products in this paper are taken over the manifold with corners $Q$.
    The maps $\pi:E\to Q$ that occur in them are mainly of two kinds, the projection of a fiber bundle and the quotient map of a torus action.
\end{remark}

\begin{theorem}\label{THM:1}
    Let $Q$ be a nice manifold with corners with facets $F_1, \ldots, F_m$, and let $\zeta\to Q$ be a principal $\T^n$-bundle with characteristic class $c(\zeta)\in H^2(Q;\Z^n)$. Let $\T^{n+m}=\T^m\times\T^n$ act on $\mathscr L_{(Q,\zeta)}$ through the induced action on the $\T^m$-coordinates and the principal action on $\zeta$. There is a graded ring isomorphism
    \begin{equation}
        H^*\bigl((\mathscr L_{(Q, \zeta)})_{\T^{n+m}};\Z\bigr)\cong\Z\langle Q\rangle.
    \end{equation}
    Write $H^*(B\T^{n+m};\Z)=\Z[x_1, \ldots, x_m, y_1, \ldots, y_n]$, with all the generators having degree $2$. Under this isomorphism, the module structure is given by $\star$-multiplication through the graded ring homomorphism
    \begin{equation}\label{eq:intro-module-one}
        x_i\longmapsto 1\otimes x_i,
        \qquad y_k\longmapsto c(\zeta)_k\otimes1.
    \end{equation}
\end{theorem}

Every locally standard $\T^n$-manifold $M$ with orbit space $Q$ is equivariantly diffeomorphic to $M_{(Q, \zeta, \lambda)}$ (see \Cref{def:quotientconstruction}) for a suitable principal $\T^n$-bundle $\zeta$ over $Q$. See \cite[Construction~8.3 and Proposition~7.8]{karshon2025classification}. 
In \Cref{sec:torus-manifolds} we show that $M_{(Q, \zeta, \lambda)}$ can be realized as a free torus quotient of the generalized moment-angle manifold $\mathscr L_{(Q,\zeta)}=\zeta\times_Q\mathscr L_Q$ of \Cref{THM:1}.
Here $\zeta$ is the principal $\T^n$-bundle determined by the free orbits of the $\T^n$-action on $M$.
This identification is what allows us to generalize the results about the equivariant cohomology of generalized moment-angle manifolds to general locally standard $\T^n$-manifolds. 
It is also the reason why we state \Cref{THM:1} in terms of the fiber product $\zeta\times_Q\mathscr L_Q$. 

Our second main result computes the equivariant cohomology ring of $M$ and its module structure (\Cref{THM:2}). We first fix the notation needed to state it precisely.
Write $Q^\circ$ for the interior of $Q$ and $M^\circ$ for its inverse image in $M$, which consists of the free orbits of $M$ under the torus action. For each facet $F_j$, choose a primitive vector $\lambda(F_j)\in\Z^n$ corresponding to the isotropy group of its interior. This is the characteristic function of $M$, and we write $\lambda(F_j)=(\lambda_{1j}, \ldots, \lambda_{nj})$.

\begin{theorem}\label{THM:2}
    Let $M$ be a smooth locally standard $\T^n$-manifold with characteristic function $\lambda$ and orbit space $Q$, where $Q$ is a nice manifold with corners with facets $F_1, \ldots, F_m$. 
    Let $\zeta\to Q$ be a principal $\T^n$-bundle restricting to the principal orbit bundle $M^\circ\to Q^\circ$. Its characteristic class
    \[
        c(\zeta)\in H^2(Q;\Z^n)
    \]
    corresponds to that of $M^\circ\to Q^\circ$ under the restriction isomorphism $H^2(Q;\Z^n)\cong H^2(Q^\circ;\Z^n)$. There is a graded ring isomorphism
    \begin{equation}
        H^*(M_{\T^n};\Z)\cong\Z\langle Q\rangle.
    \end{equation}
    For the standard degree-two generators $y_1, \ldots, y_n$ of $H^*(B\T^n;\Z)$, the module structure is given by $\star$-multiplication through
    \begin{equation}\label{eq:intro-module-two}
        y_i\longmapsto c(\zeta)_i\otimes1+\sum_{j=1}^m\lambda_{ij}(1\otimes x_j).
    \end{equation}
\end{theorem}

The paper is organized as follows.
\Cref{sec:torus-manifolds} relates generalized moment-angle manifolds to locally standard $\T^n$-manifolds and deduces \Cref{THM:2} from \Cref{THM:1}.
\Cref{sec:preliminary} collects the basic notions and properties of homotopy colimits.
In \Cref{sec:gppomwc} we generalize the notion of the polyhedral product. We show that every generalized moment-angle manifold is equivariantly homeomorphic to such a product over a nice manifold with corners (\Cref{thm:mamehtpp}).
\Cref{sec:homotopycolimit} identifies the generalized polyhedral product as a homotopy colimit and reduces its equivariant cohomology to that of a simpler diagram.
\Cref{thm:generalpolyhedral} states that the equivariant cohomology ring of this product is an inverse limit, and describes its module structure in terms of the characteristic classes.
\Cref{sec:computation} applies the Segal spectral sequence to the homotopy colimit of \Cref{sec:homotopycolimit}. The $E_2$-page is given by the derived limits of the diagram of cohomology groups. 
We show that the higher derived limits vanish. \Cref{thm:generalpolyhedral} then follows from the collapsing of the spectral sequence, and \Cref{THM:1} follows by identifying the polyhedral product with $\mathscr L_{(Q, \zeta)}$ as in \Cref{thm:mamehtpp}. 
The quasi-cubical decomposition used in \Cref{sec:gppomwc} is constructed in the appendix.

    \section{From Theorem 1 to Theorem 2}\label{sec:torus-manifolds}
In this section, we prove \Cref{THM:2} by applying \Cref{THM:1} to a fiber product with the principal orbit bundle and taking a free torus quotient.

Let $M$ be a smooth locally standard $\T^n$-manifold with orbit space $Q$ and characteristic function $\lambda:\mathcal{F}\to\Z^n$. The principal orbit bundle $M^\circ\to Q^\circ$ extends to a smooth principal $\T^n$-bundle $\zeta\to Q$, and there is a $\T^n$-equivariant diffeomorphism
\[
    M\cong M_{(Q, \zeta, \lambda)}
\]
over the identity of $Q$. See \cite[Construction~8.3 and Proposition~7.8]{karshon2025classification}. The equivalence relation in \cite[Construction~8.3]{karshon2025classification} is exactly \eqref{eq:relationdefquotient}. Write $c(\zeta)\in H^2(Q;\Z^n)$ for its characteristic class, whose pullback to $H^2(Q^\circ;\Z^n)$ is the class of the principal bundle $M^\circ\to Q^\circ$ \cite[Lemma~3.13]{karshon2025classification}.

Recall that $\mathscr L_{(Q,\zeta)}=\mathscr L_Q\times_Q\zeta$, on which $\T^{n+m}=\T^m\times\T^n$ acts through the two factors as in \Cref{THM:1}. Define
\begin{equation*}
    \Phi:\T^{n+m}\longrightarrow\T^n, \quad (t, s)\longmapsto\Lambda(t)s.
\end{equation*} The homomorphism $\Phi$ has a section $s\mapsto(1, s)$, and
\begin{equation*}
    K:=\ker\Phi=\{(t, \Lambda(t)^{-1}):t\in\T^m\}\cong\T^m.
\end{equation*}
This construction does not require a splitting of $\lambda$.

\begin{lemma}\label{LEM:EquivariantMap}
    The $K$-action on $\mathscr L_{(Q,\zeta)}$ is free, and there is a $\T^{n+m}/K\cong\T^n$-equivariant homeomorphism
    \begin{equation}\label{eq:quotientismomentangle}
        \mathscr L_{(Q,\zeta)}/K\cong M_{(Q, \zeta, \lambda)}.
    \end{equation}
\end{lemma}

\begin{proof}
    For $q\in Q$, let $J(q)=\{j:q\in F_j\}$.
    By local standardness, $\{\lambda(F_j):j\in J(q)\}$ is part of an integral basis of $\Z^n$, so $\Lambda|_{\T^m_{J(q)}}$ is injective.
    Thus the isotropy subgroup $\T^m_{J(q)}\times\{1\}$ at any point over $q$ meets $K$ trivially, and $K$ acts freely.

    By definition, $\mathscr L_{(Q,\zeta)}$ and $M_{(Q, \zeta, \lambda)}$ are the quotients of $\T^m\times\zeta$ and $\zeta$ by the respective relations
    \[
        \begin{aligned}
            (t, p)&\sim(ht, p), && h\in\T^m_{J(q)}, \\
            p&\sim gp, && g\in\T^n_{\lambda, J(q)},
        \end{aligned}
        \qquad q=\pi_\zeta(p).
    \]
    We write $[t, p]$ and $[p]$ for the corresponding equivalence classes.
    Since $\Lambda(\T^m_{J(q)})=\T^n_{\lambda, J(q)}$, the map
    \begin{equation}\label{eq:quotientmap}
        \mathscr L_{(Q,\zeta)}\longrightarrow M_{(Q, \zeta, \lambda)}, \qquad
        [t, p]\longmapsto[\Lambda(t)p]
    \end{equation}
    is well-defined. It is $\Phi$-equivariant and hence descends to $\mathscr L_{(Q,\zeta)}/K\to M_{(Q, \zeta, \lambda)}$.

    The induced map has inverse
    \[
        M_{(Q, \zeta, \lambda)}\longrightarrow \mathscr L_{(Q,\zeta)}/K, \qquad
        [p]\longmapsto K\cdot[1, p].
    \]
    This is well-defined, since replacing $p$ by $\Lambda(h)p$ with $h\in\T^m_{J(q)}$ leaves the $K$-orbit unchanged:
    \[
        (h^{-1}, \Lambda(h))\cdot[1, p]
        =[h^{-1}, \Lambda(h)p]=[1, \Lambda(h)p].
    \]
    The two maps are continuous by the quotient topology and mutually inverse, giving the asserted equivariant homeomorphism.
\end{proof}

We next compare the Borel constructions, keeping track of their maps to the classifying spaces to determine the module structure.
\begin{corollary}\label{cor:diagram}
    There is a commutative diagram
    \begin{equation}\label{eq:introcommdiag}
    \begin{tikzcd}
        {(\mathscr L_{(Q,\zeta)})_{\T^{n+m}}} & {(M_{(Q,\zeta,\lambda)})_{\T^n}} \\
        {B\T^{n+m}} & {B\T^n}
        \arrow["\simeq", from=1-1, to=1-2]
        \arrow[from=1-1, to=2-1]
        \arrow[from=1-2, to=2-2]
        \arrow["{B\Phi}"', from=2-1, to=2-2]
    \end{tikzcd}
    \end{equation}
    whose top arrow is a homotopy equivalence.
\end{corollary}

\begin{proof}
    The splitting of $\Phi$ identifies $\T^{n+m}$ with $K\times\T^n$. Choose $E\T^{n+m}=EK\times E\T^n$ with respect to this identification. Then the top arrow is the natural map
    \begin{equation*}
        E\T^n\times_{\T^n}(EK\times_K \mathscr L_{(Q,\zeta)})\longrightarrow E\T^n\times_{\T^n}(\mathscr L_{(Q,\zeta)}/K),
    \end{equation*}
    followed by the homeomorphism from \Cref{LEM:EquivariantMap}. With these identifications, $B\Phi$ is projection onto the second factor, so the diagram commutes.

    Since the action of $K$ on $\mathscr L_{(Q,\zeta)}$ is free, $EK\times_K \mathscr L_{(Q,\zeta)} \to \mathscr L_{(Q,\zeta)}/K$ is a bundle with fiber $EK$. 
    Hence the top arrow of \eqref{eq:introcommdiag} is a bundle with contractible fiber $EK$, and thus a homotopy equivalence.
\end{proof}

\begin{remark}\label{rem:relationtoYu}
    In the fiber product $\mathscr L_{(Q,\zeta)}=\zeta\times_Q\mathscr L_Q$, the group $\T^n$ acts freely on $\zeta$, and the quotient is $\mathscr L_Q$. Since the Borel construction of the principal $\T^n$-bundle $\zeta$ is homotopy equivalent to $Q$, this gives a homotopy equivalence
    \[
        (\mathscr L_{(Q,\zeta)})_{\T^{n+m}}\simeq(\mathscr L_Q)_{\T^m}.
    \]
    Together with \Cref{LEM:EquivariantMap} and \Cref{cor:diagram}, we obtain
    \[
        (M_{(Q, \zeta, \lambda)})_{\T^n}\simeq(\mathscr L_Q)_{\T^m}
    \]
    for every locally standard $\T^n$-manifold $M$. The equivariant cohomology ring of $\mathscr L_Q$ has been computed by Yu in \cite{MR4721373}, so as far as the ring structure is concerned, \Cref{THM:2} already follows. The module structure, however, cannot be obtained directly from these homotopy equivalences.
\end{remark}

We now use this comparison to deduce \Cref{THM:2} from \Cref{THM:1}.
\begin{proof}[Proof of \Cref{THM:2}]
    The equivariant homeomorphism $M\cong M_{(Q, \zeta, \lambda)}$ induces a homeomorphism of Borel constructions over $B\T^n$.
    Applying \Cref{cor:diagram} and \Cref{THM:1}, we obtain graded ring isomorphisms
    \begin{equation*}
        H^*((M_{(Q, \zeta, \lambda)})_{\T^n};\Z)
        \cong H^*((\mathscr L_{(Q,\zeta)})_{\T^{n+m}};\Z)
        \cong\Z\langle Q\rangle.
    \end{equation*}

    It remains to determine the module structure, for which we compute the pullback $(B\Phi)^*$.
    The homomorphism $\Phi$ induces the matrix $(\lambda\mid I_n)$ on fundamental groups. By the natural isomorphisms $\pi_2(B\T^r)\cong\pi_1(\T^r)$, the Hurewicz theorem and the universal coefficient theorem, its pullback on degree-two cohomology is
    \begin{equation}\label{eq:pullbackofPhi}
        (B\Phi)^*(y_i)=\sum_{j=1}^m\lambda_{ij}x_j+y_i\in H^2(B\T^{n+m};\Z),
    \end{equation}
    where the generators on the right are those of \Cref{THM:1}. Thus \eqref{eq:intro-module-one} sends the right-hand side in \eqref{eq:pullbackofPhi} to $c(\zeta)_i\otimes1+\sum_{j=1}^m\lambda_{ij}(1\otimes x_j)$. By the commutativity of \eqref{eq:introcommdiag}, this is the image of $y_i$ under the module homomorphism, proving \eqref{eq:intro-module-two}. Explicitly, for $\alpha\otimes f\in H^*(F_{\cap J};\Z)\otimes\Z[\mathbf{x}_J^+]$, its action is
    \begin{equation}\label{eq:expanded-module-two}
        y_i\cdot(\alpha\otimes f)=
        (c(\zeta)_i|_{F_{\cap J}}\cup\alpha)\otimes f
        +\sum_{j=1}^m\lambda_{ij}\alpha|_{F_{\cap(J\cup\{j\})}}\otimes x_jf.
    \end{equation}
\end{proof}

    \section{Preliminary}\label{sec:preliminary}
We introduce some necessary concepts about \emph{homotopy colimits} in this section. The material is mainly taken from \cite{dugger2008primer,MR1944041,MR1711612}.  
A detailed discussion of homotopy colimits requires some knowledge of \emph{model category theory} and \emph{simplicial homotopy theory}. Because these theories are quite technical, we do not attempt to cover them here. The statements in this section are therefore presented without proof, but we give enough references for the reader to verify them.
We recommend the reader who is not familiar with homotopy colimits to read \cite{dugger2008primer} to get a quick start.

\begin{definition}
    Let $\mathcal{M}$ be a category. 
For a small category $I$ called the \emph{indexing category}, an object of the functor category $\mathcal{M}^I$ is a functor $D: I\to \mathcal{M}$ called a \emph{diagram} in $\mathcal{M}$ of type $I$. 
We call the image $D(i\to j)$ of a morphism $i\to j$ of $I$ a \emph{structure map} of the diagram $D$. 
A \emph{constant diagram} is a diagram $D$ such that: for each $i\in I$, $D(i)$ is a fixed object $W\in \mathcal{M}$, and for each $i\to j$, $D(i\to j)$ is the identity morphism. We denote such a constant diagram by $c(W)$. 
A map between two diagrams of the same type is a natural transformation. 
\end{definition}

The definition of the homotopy colimit of a diagram requires the notion of the \emph{(co)simplicial object} and its \emph{geometric realization}. We introduce these concepts briefly in the following and refer the reader to \cite[pp. 3~9]{MR1711612} for details. 

Simplicial objects and cosimplicial objects are specific diagrams in a given category. Their definitions are dual to each other. Below we introduce the cosimplicial object using standard geometric simplices.
\begin{example}
    Let $\Delta^n$ be the standard geometric $n$-simplex. Let $d^i: \Delta^n \to \Delta^{n+1}$ be the map which identifies $\Delta^n$ as the $i$-th face of $\Delta^{n+1}$ and 
    $s^i: \Delta^{n+1} \to \Delta^n$ be the affine map that replaces the adjacent barycentric coordinates $t_i, t_{i+1}$ by their sum $t_i+t_{i+1}$. We call $d^i$ a \emph{coface map} and $s^i$ a \emph{codegeneracy map}. These maps satisfy the so-called \emph{cosimplicial identities} (see \cite[p. 4]{MR1711612}). The collection of all standard geometric $n$-simplices together with \emph{coface maps} and \emph{codegeneracy maps} forms a cosimplicial object in $\mathbf{Top}$.
\end{example}
A simplicial object in $\mathbf{Top}$, also known as a \emph{simplicial space}, is defined dually to the above example. Roughly speaking, a simplicial space is a sequence of spaces $X_n$ labeled by $n\in \mathbb{N}$. There is a set of continuous maps $d_i:X_n\to X_{n-1}$ ($0\leq i \leq n$) called \emph{face maps} and $s_i: X_{n}\to X_{n+1}$ ($0 \leq i\leq n$) called \emph{degeneracy maps}. These maps satisfy the \emph{simplicial identities} (see \cite[p. 5]{MR1711612}).

Given a diagram $D: I\to \mathbf{Top}$, the \emph{simplicial replacement} of $D$ is a \emph{simplicial space} $serp(D)$ defined by
\begin{equation*}
    serp(D)_n = \coprod_{i_0 \to i_1 \to \cdots \to i_n} D(i_0).
\end{equation*}
Here, the notation \(\coprod_{i_0 \to i_1 \to \cdots \to i_n}\) denotes the disjoint union over all chains of morphisms in $I$ of length \(n\).
For the definition of face maps $d_i$ and degeneracy maps $s_i$ of this simplicial space, the reader can refer to \cite[p. 16]{dugger2008primer}. 

\begin{definition}\label{def:homotopycolimit}
    The \emph{homotopy colimit} of the diagram $D: I \to \mathbf{Top}$ is defined by 
    \begin{equation}\label{def:hocolimdef}
        \begin{gathered}
            \hocolim_I D:=\left(\coprod_{n\in \mathbb{N}} serp(D)_n \times \Delta^n\right)\Big/\sim, \\ 
            (d_i(x), p) \sim (x, d^i(p)), \qquad (s_i(x), p) \sim (x, s^i(p)). 
        \end{gathered}
    \end{equation}
    This definition comes from a general construction of the \emph{realization of a simplicial object}. 
\end{definition}

Here are observations that follow directly from the definition of the homotopy colimit, which we will use later:
\begin{enumerate}
    \item If the diagram $D$ admits the action of a group $G$, that is, the spaces in the diagram are $G$-spaces and the maps are $G$-equivariant, then $\hocolim_I D$ admits a $G$-action: one checks that $serp(D)$ admits a $G$-action with equivariant face and degeneracy maps, and equipping each $\Delta^n$ with the trivial $G$-action realizes $\hocolim_I D$ as a $G$-space.
    \item There is a natural filtration $F_0\subset F_1 \subset \cdots F_k \subset \cdots$ of the homotopy colimit $\hocolim_I D$. Here, $F_k$ is defined to be the image of 
    \begin{equation*}
        \coprod_{0\leq n \leq k} serp(D)_n \times \Delta^n
    \end{equation*}
    in $\hocolim_I D$ under the quotient map. Note that $F_0$ is the disjoint union of spaces in the diagram. 
    A spectral sequence can be constructed from this filtration, which converges to the cohomology of $\hocolim_I D$. This spectral sequence was introduced by Segal in \cite{MR232393}. 
     We shall use this spectral sequence to carry out computations in \Cref{sec:computation}. 
    \item Let $W$ be a topological space. If there is a natural transformation from $D$ to the constant diagram $c(W)$, then there exist induced maps $serp(D)_n \to W$ which commute with face and degeneracy maps. By considering the composition $serp(D)_n \times \Delta^n \to serp(D)_n \to W$ for each $n$, we get a map $\hocolim_I D \to W$.
\end{enumerate}

We use weak equivalences in the sense of weak homotopy equivalences, and we say that a space is cofibrant when it is a retract of a CW-complex. See \cite[Definition~7.1.6]{MR1944041}. A CW-complex is cofibrant, and the inclusion of a CW-subcomplex is a cofibration.

The following lemma explains our motivation for introducing the homotopy colimit: we wish to work with homotopies more freely, without changing the homotopy type.
\begin{lemma}[{\cite[Theorem~18.5.3(1)]{MR1944041}}]\label{prop:lemmahomotopycolimit}
    Let $X, Y$ be diagrams in $\mathbf{Top}$ of type $I$. Then the map $\hocolim f: \hocolim X \to \hocolim Y$ induced by $f: X\to Y$ is a weak equivalence if $f_i:X_i\to Y_i$ is a weak equivalence between cofibrant objects for each $i\in I$. 
\end{lemma}
A natural transformation is called a \emph{pointwise weak equivalence} if each of its component maps is a weak equivalence. The above lemma says that a pointwise weak equivalence between diagrams of cofibrant spaces induces a weak equivalence of homotopy colimits. For colimits, however, this is generally not true.

According to \cite[p. 11]{dugger2008primer}, if we collapse the simplices in the homotopy colimit $\hocolim_IX$ (see \Cref{def:hocolimdef}), then we get a collapsing map $\hocolim_I X \to \colim_I X$.
In general, this collapsing map is not a weak equivalence. When the index set $I$ is induced by a poset, there is a criterion for whether this map is a weak equivalence.

A poset $(I, \leq)$ determines a category whose objects are the elements of $I$, and there is a morphism from $i$ to $j$ if and only if $i \leq j$. We identify $I$ with $(I, \leq)$. For a poset $I$ and an element $i\in I$, we use $I_{<i}$ to denote the full subcategory of $I$ with objects $\{j \in I | j<i\}$. 
Let $X$ be a diagram of type $I$ where $I$ is a poset. 
Then for each $j\in I_{<i}$, there is a map $X_j \to X_i$. 
The universal property of the colimit functor gives a map called a \emph{latching map}:
\begin{equation}
    \colim_{I_{<i}} X \to X_i. 
\end{equation}
The following lemma provides a sufficient condition for the collapsing map to be a weak equivalence.
\begin{lemma}[{\cite[Theorem~19.9.1]{MR1944041}}]\label{lemma:latchinglemma}
    Let $I$ be a finite poset and $X$ be a diagram in $\mathbf{Top}$ of type $I$. 
    If for each $i\in I$, the latching map $\colim_{I_{<i}} X \to X_i$ is a cofibration, then the collapsing map $\hocolim_I X \to \colim_I X$ is a weak equivalence. 
\end{lemma}

    \section{Generalized moment-angle manifolds as polyhedral products}\label{sec:gppomwc}
Polyhedral products are used in \cite[Chapter~4]{MR3363157} to study the cohomology and equivariant cohomology rings of moment-angle complexes and manifolds. Yu also used this approach in \cite{MR4721373}, where the principal torus bundle over a nice manifold with corners is trivial. We follow a similar approach here.

\subsection{Quasi-cubical decompositions and rim-cubicalization}\label{section:cubicalization}
\begin{figure}[htbp]                
    \centering                      
    \includegraphics[width=0.75\textwidth]{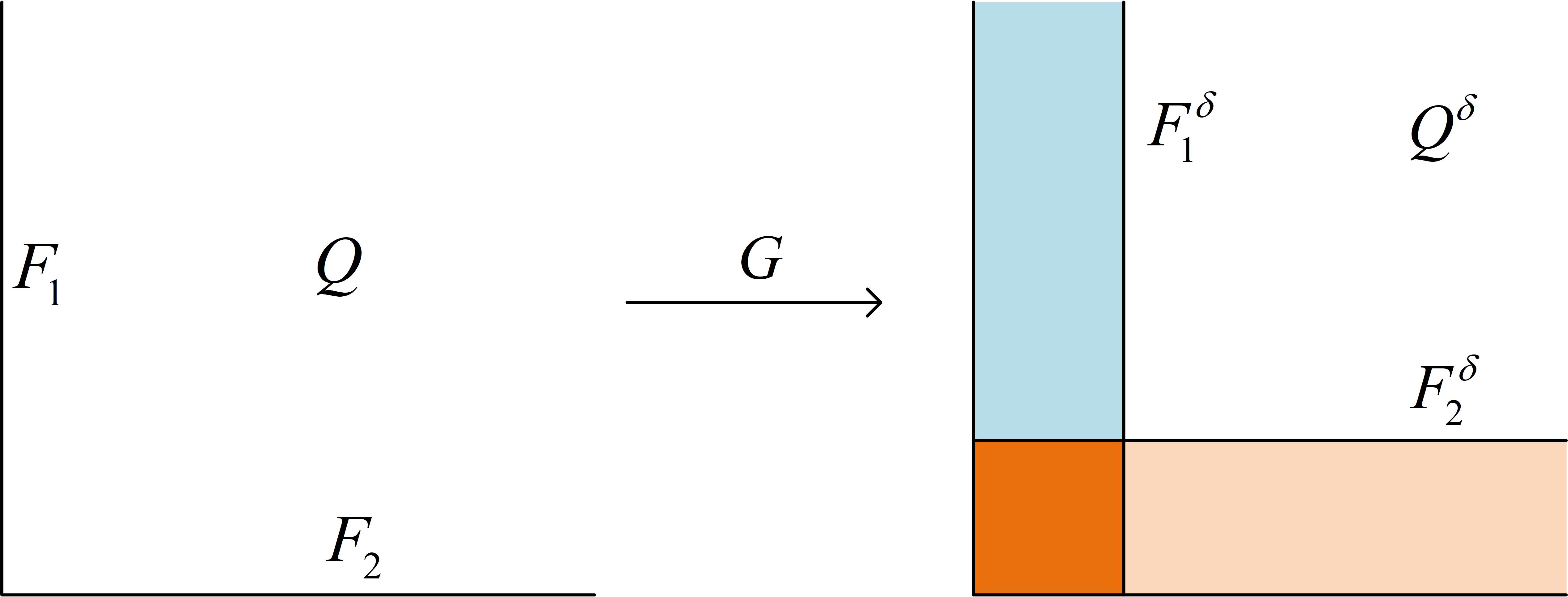} 
    \caption{The quasi-cubical decomposition of $Q$ with two facets}               
    \label{fig:quasi-cubical decomposition}            
\end{figure}

We begin with a decomposition of $Q$ that is analogous to the collar neighborhood theorem. By the collar neighborhood theorem, a manifold $W$ with boundary can be decomposed into $\partial W\times[0, \delta]$ and its complement $W\setminus(\partial W\times[0, \delta))$ for some $\delta>0$. A \emph{quasi-cubical decomposition} of a nice manifold with corners is the analogue of this decomposition in which the boundary is replaced by all the faces. The number $\delta$ is the width of the collars, and $Q^{\delta}$ denotes the copy of $Q$ obtained by pushing $Q$ away from its faces by this amount. 
Moreover, $Q$ is written as a union of blocks indexed by the subsets $J\subseteq[m]$, the block indexed by $J$ being the product $F^{\delta}_{\cap J} \times [0, \delta]^{|J|}$, where $|J|$ denotes the cardinality of $J$. 
Here, $F^{\delta}_{\cap J}$ is homeomorphic to $F_{\cap J}$ for each $J\subseteq [m]$. Intuitively, the blocks are drawn in \Cref{fig:quasi-cubical decomposition}. 

The difference from the collar neighborhood theorem is that the collars of different facets have to be compatible with one another, since they meet near the faces of higher codimension. This makes the construction delicate and the proof rather technical, so we give the precise definition and the existence proof in Appendix~\ref{app:quasicubical}.
Yu used a decomposition of the same type in \cite{MR4721373}. However, the proof of \cite[Lemma~2.1]{MR4721373} is not correct.

\Cref{lemma:embeddingofQandxi} below is the corresponding statement for our setting, and its proof uses the decomposition of Appendix~\ref{app:quasicubical}. In what follows we take $\delta=1$ for convenience.
Given a quasi-cubical decomposition, one embeds $Q$ into $Q\times[0, 1]^m$ as follows.
\begin{lemma}\label{lemma:embeddingofQandxi}
    Let $Q$ be a nice manifold with corners. Then there is an embedding $\hat{f}: Q \to Q\times [0, 1]^m$ whose image is of the form:
    \begin{equation}
         \hat{f}(Q) = \widehat{Q}:=\bigcup_{J \subseteq [m]} F_{\cap J}\times \prod_{j\in J} [0, 1]_{(j)} \times \prod_{j\notin J} 1_{(j)}.
    \end{equation}
    Moreover, for any $J\neq \varnothing$, we have 
    \begin{equation}\label{eq:preimageoff}
            \hat{f}(F_{\cap J})=\widehat{F}_{\cap J}:=\bigcup_{J \subseteq J'} F_{\cap J'} \times \prod_{j\in J} 0_{(j)} \times \prod_{j \in J'\setminus J} [0, 1]_{(j)} \times \prod_{j\notin J'} 1_{(j)}.
    \end{equation}
    Furthermore, we have the following principal bundle map which is an embedding:
    \begin{equation}\label{cd:pullbackbundle}
        \begin{tikzcd}
            \xi & {\xi\times[0, 1]^m} \\
            Q & {Q\times [0, 1]^m}
            \arrow["{\hat{f}_*}", from=1-1, to=1-2]
            \arrow["{\pi_\xi}"', from=1-1, to=2-1]
            \arrow["{\pi_{\xi} \times \id}", from=1-2, to=2-2]
            \arrow["\hat{f}"', from=2-1, to=2-2]
        \end{tikzcd}
    \end{equation}
\end{lemma}
We call $\widehat{Q}$ the rim-cubicalization of $Q$ in $Q\times [0, 1]^m$.

\subsection{The generalized polyhedral product}
In this subsection we introduce the generalized polyhedral product, which generalizes both the polyhedral product in \cite[Section~4.2.1]{MR3363157} and the construction of Yu in \cite[Section~2.2]{MR4721373}. We then show in \Cref{thm:mamehtpp} that the generalized moment-angle manifold is equivariantly homeomorphic to the generalized polyhedral product.

Take ${\mathrm{D}^2}\subseteq\C$ to be the closed unit disk and $\mathrm{S}^1=\partial{\mathrm{D}^2}$. The product $({\mathrm{D}^2})^m$ carries the coordinatewise action of $\T^m$
\begin{equation}\label{def:multiplicationoftonc}
    (g_1, \ldots, g_m)\cdot(z_1, \ldots, z_m)=(g_1z_1, \ldots, g_mz_m).
\end{equation}
Recall that for $J\subseteq [m]$, the space $(\mathrm{D}^2, \mathrm{S}^1)^J$ is defined by 
\begin{equation*}
    (\mathrm{D}^2, \mathrm{S}^1)^J := \prod_{j\in J}{\mathrm{D}^2_{(j)}}\times\prod_{j\notin J}\mathrm{S}^1_{(j)}
\end{equation*}

\begin{definition}\label{def:genpolyprod}
    Let $\zeta\to Q$ be a principal $\T^n$-bundle. For $J\subseteq[m]$, we set
    \begin{equation}\label{eq:defofoneplyhedra}
        ({\mathrm{D}^2}, \mathrm{S}^1)^{(\zeta, J)}
        :=\zeta|_{F_{\cap J}}\times (\mathrm{D}^2, \mathrm{S}^1)^J  =\zeta|_{F_{\cap J}}\times\prod_{j\in J}{\mathrm{D}^2_{(j)}}\times\prod_{j\notin J}\mathrm{S}^1_{(j)}
    \end{equation}
    We call $(\mathrm{D}^2, \mathrm{S}^1)^{(\zeta, J)}$ the \emph{subproduct} indexed by $J$, and define the \emph{generalized polyhedral product}
    \begin{equation}\label{eq:compbypolyhedraldef}
        ({\mathrm{D}^2}, \mathrm{S}^1)^{\zeta}:=\bigcup_{J\subseteq[m]}({\mathrm{D}^2}, \mathrm{S}^1)^{(\zeta, J)},
    \end{equation}
    which is a subspace of $\zeta \times (\mathrm{D}^2)^m$. 
\end{definition}

The group $\T^{n+m}$ acts on $\zeta \times (\mathrm{D}^2)^m$ by
\begin{equation}\label{eq:actonassociatedbundle}
    (g, h)\cdot(v, z)=\bigl(hv, gz\bigr), \qquad (g, h)\in\T^m\times\T^n, \quad (v, z) \in \zeta \times (\mathrm{D}^2)^m,
\end{equation}
where $hv$ denotes the principal action of $\T^n$ on $\zeta$, and $gz$ denotes the action defined in \eqref{def:multiplicationoftonc}. Each subproduct is an invariant subspace of $\zeta \times (\mathrm{D}^2)^m$, and so is the generalized polyhedral product. 

Hence the subproducts and the generalized polyhedral product inherit the action of $\T^{n+m}$ on $\zeta \times (\mathrm{D}^2)^m$.
The orbit space of $\mathrm{D}^2$ under the action of $\T^1$ is the interval $[0, 1]$. Under this quotient map, $\mathrm{S}^1$ is mapped to $1 \in [0, 1]$. 
Thus, the orbit space of $\zeta \times (\mathrm{D}^2)^m$ under the action of $\T^{n+m}$ is $Q\times [0, 1]^m$,
and for each $J\subseteq [m]$, the orbit space of the subproduct indexed by $J$ is 
\[
    F_{\cap J} \times \prod_{j\in J} [0, 1]_{(j)} \times \prod_{j\notin J} 1_{(j)}.
\]
Hence, the orbit space of $(\mathrm{D}^2, \mathrm{S}^1)^{\zeta}$ under the action of $\T^{n+m}$ is
\[
    \widehat{Q} = \bigcup_{J \subseteq [m]} F_{\cap J}\times \prod_{j\in J} [0, 1]_{(j)} \times \prod_{j\notin J} 1_{(j)}.
\]

To prepare for the proof of \Cref{thm:mamehtpp}, we describe $(\mathrm{D}^2)^m$, and its subspaces $(\mathrm{D}^2, \mathrm{S}^1)^J$, as quotients of $\T^m\times[0,1]^m$. In $\C$, the interval $[0, 1]$ is the intersection of $\mathrm{D}^2$ with the real axis. Take $[0, 1]^m$ as a subspace of $(\mathrm{D}^2)^m$. Then we can define a map $\nu$ by the multiplication given in 
\begin{equation}
    \begin{aligned}
        \nu: \T^m \times [0, 1]^m &\to (\mathrm{D}^2)^m \\ 
        \bigl((g_1, \cdots, g_m), (r_1, \cdots, r_m)\bigr) &\mapsto (g_1r_1, \cdots, g_mr_m). 
    \end{aligned}
\end{equation}
For $r=(r_1, \cdots, r_m)\in [0, 1]^m$, define 
\begin{equation}
    Z(r):= \{i\in[m] \mid r_i =0 \}. 
\end{equation}
Then for $r, r' \in [0, 1]^m$, there exist $g, g' \in \T^m$ such that $\nu(g, r) = \nu(g', r')$ if and only if 
\begin{equation}\label{eq:quotientdefofmu}
    r=r' \text{ and } g(g')^{-1} \in \T^m_{Z(r)},
\end{equation}
Note that $\nu$ is surjective, hence we can take $(\mathrm{D}^2)^m$ as the quotient space of 
$\T^m \times [0, 1]^m$ under the relation in \eqref{eq:quotientdefofmu}. For each $J\subseteq [m]$, the image of the subspace
\begin{equation}
    \T^m \times \prod_{j\in J} [0, 1]_{(j)} \times \prod_{j\notin J} 1_{(j)} 
    = \prod_{j\in J} \bigl(\T^1_{(j)} \times [0, 1]_{(j)}\bigr) \times \prod_{j\notin J}\bigl(\T^1_{(j)} \times  1_{(j)}\bigr)
\end{equation}
under $\nu$ is exactly $(\mathrm{D}^2, \mathrm{S}^1)^J$. 

\begin{theorem}\label{thm:mamehtpp}
    There is a $\T^{n+m}$-equivariant homeomorphism
    \begin{equation}\label{eq:zethomeomorphism}
        \scrl_{(Q, \zeta)}\cong\;({\mathrm{D}^2}, \mathrm{S}^1)^{\zeta}.
    \end{equation}
\end{theorem}

\begin{proof}
    Setting $\xi = \zeta \times \T^m$ in \Cref{lemma:embeddingofQandxi}, we obtain the following commutative diagram:
    \[\begin{tikzcd}
        {\zeta \times \T^m} & {\zeta \times \T^m \times [0, 1]^m} \\
        Q & {Q\times [0, 1]^m}
        \arrow["{\hat{f}_*}", from=1-1, to=1-2]
        \arrow[from=1-1, to=2-1]
        \arrow[from=1-2, to=2-2]
        \arrow["{\hat{f}}", from=2-1, to=2-2]
    \end{tikzcd}.\]
    By the pullback property, we may choose $\hat{f}_*$ of the form
    \begin{equation}
        \hat{f}_*(v, g) = (\phi(v), g, T(v)),
    \end{equation}
    where $\phi:\zeta\to\zeta$ is a $\T^n$-equivariant bundle map covering $p_1\circ\hat{f}$ and $T: \zeta \to [0, 1]^m$ is the composition:
    \begin{equation}\label{eq:imageoffstar}
        \zeta \stackrel{\pi_{\zeta}}{\longrightarrow} Q \stackrel{\hat{f}}{\longrightarrow} Q \times [0, 1]^m \stackrel{p_2}{\longrightarrow} [0, 1]^m. 
    \end{equation}
    By the definition of the generalized moment-angle manifold $\scrl_{(Q, \zeta)}$, it is the quotient space $(\zeta \times \T^m)/ \sim_1$, where $J$ denotes the unique subset with $\pi_{\zeta}(v)\in F_{\cap J}^{\circ}$ and
    \begin{equation}\label{eq:quotientconstructionofmomentangle}
        (v, g) \sim_1 (v', g') \Leftrightarrow v=v' \text{ and } g(g')^{-1} \in \T^m_{J}.
    \end{equation}
    Since $F_{\cap J}^{\circ}=F_{\cap J}\setminus\bigcup_{J\subsetneq J'}F_{\cap J'}$ and $\hat{f}$ is an embedding, \eqref{eq:preimageoff} gives
    \begin{equation*}
        \hat{f}(F^{\circ}_{\cap J}) = \hat{f}(F_{\cap J})\setminus\bigcup_{J\subsetneq J'}\hat{f}(F_{\cap J'}) = \{(q,r)\in\widehat Q \mid Z(r)=J\}.
    \end{equation*}
    Hence, for any $v \in \zeta$, we have 
    \begin{equation*}
        Z(T(v)) = J \Leftrightarrow \pi_{\zeta}(v) \in F_{\cap J}^{\circ}. 
    \end{equation*}
    By the embedding $\hat{f}_*:\zeta\times\T^m\to\zeta\times\T^m\times[0,1]^m$, the relation given in \eqref{eq:quotientconstructionofmomentangle} can be rewritten as a relation on the subspace $\hat{f}_*(\zeta\times\T^m)$:
    \begin{equation*}
        (\phi(v), g, T(v)) \sim_2 (\phi(v'), g', T(v')) \Leftrightarrow v=v' \text{ and } g(g')^{-1} \in \T^m_{Z(T(v))}.
    \end{equation*}
    Comparing the relation $\sim_2$ with \eqref{eq:quotientdefofmu}, two points of $\hat{f}_*(\zeta \times \T^m)$ have the same image under $\id \times \nu$ if and only if they are related by $\sim_2$. Since $\nu$ is a quotient map, so is $\id \times \nu$, and hence it induces a homeomorphism
    \[
        \hat{f}_*(\zeta \times \T^m)/\sim_2 \;\cong\; (\id \times \nu)\bigl(\hat{f}_*(\zeta \times \T^m)\bigr).
    \]
    By the description of $(\mathrm{D}^2, \mathrm{S}^1)^J$ above, $(\id \times \nu)\bigl(\hat{f}_*(\zeta \times \T^m)\bigr)=\bigcup_{J\subseteq[m]}\zeta|_{F_{\cap J}}\times(\mathrm{D}^2, \mathrm{S}^1)^J=(\mathrm{D}^2, \mathrm{S}^1)^{\zeta}$.
    Since $\hat{f}_*$ is a homeomorphism onto its image and $\sim_2$ is obtained from $\sim_1$ via $\hat{f}_*$, it induces a homeomorphism $(\zeta\times\T^m)/\sim_1\cong\hat{f}_*(\zeta\times\T^m)/\sim_2$.
Hence, we have 
\begin{equation*}
    \begin{aligned}
        \scrl_{(Q, \zeta)} &= (\zeta \times \T^m)/ \sim_1 \ \cong \hat{f}_* (\zeta \times \T^m) / \sim_2  \\
    &\cong ((\id \times \nu) \circ \hat{f}_*) (\zeta \times \T^m) = (\mathrm{D}^2, \mathrm{S}^1)^{\zeta}.
    \end{aligned}
\end{equation*}
It is straightforward to verify that each homeomorphism given above is $\T^{n+m}$-equivariant. 
\end{proof}

    \section{From the polyhedral product to a diagram of cohomology rings}\label{sec:homotopycolimit}
In this section we express the Borel construction of the generalized polyhedral product as the homotopy colimit of a diagram of spaces, and compute the corresponding diagram of cohomology rings together with its module structure.

\subsection{Refining the polyhedral product}\label{subsec:generalpolyhedral}
For the classical polyhedral product over a simplicial complex $\mathcal{K}$ on $[m]$, if $J$ and $J'$ are simplices of $\mathcal{K}$ with $J\subseteq J'$, then $(\mathrm{D}^2, \mathrm{S}^1)^{J}$ is contained in $(\mathrm{D}^2, \mathrm{S}^1)^{J'}$, so we can take the polyhedral product as the colimit of the diagram formed by these spaces and their natural inclusions (see \cite[Section~4.2.1]{MR3363157}). In our situation, however, the subproduct indexed by $J$ is $\zeta|_{F_{\cap J}}\times(\mathrm{D}^2, \mathrm{S}^1)^{J}$, and $\zeta|_{F_{\cap J}}$ shrinks as $J$ grows. Hence no inclusion analogous to the above holds between the subproducts indexed by $J$ and by $J'$ when $J\subseteq J'$. If we still want to regard $(\mathrm{D}^2, \mathrm{S}^1)^{\zeta}$ as the colimit of a diagram, we need a finer decomposition of it.   
 
\begin{definition}
    For $L\subseteq J\subseteq[m]$, define
\begin{equation}\label{eq:Xiprod}
    ({\mathrm{D}^2}, \mathrm{S}^1)^{(\zeta, J, L)}:=\zeta|_{F_{\cap J}}\times (\mathrm{D}^2, \mathrm{S}^1)^L = \zeta|_{F_{\cap J}}\times\prod_{j\in L}{\mathrm{D}^2_{(j)}}\times\prod_{j\notin L}\mathrm{S}^1_{(j)}.
\end{equation}
We call these spaces the \emph{refined subproducts} indexed by $(J, L)$. 
\end{definition}

Define a poset by
\begin{equation}\label{eq:defofP}
    \mathscr{P}:=\{ (J, L) \mid  L\subseteq J \subseteq[m]\}, 
\end{equation}
where $(J, L) \leq (J', L')$ if and only if
\[
    L\subseteq L' \subseteq J' \subseteq J   . 
\]
\begin{remark}
    This index set is similar to the one in \cite[(15)]{MR4721373}, where a pair consists of a face $f$ of $Q$ together with a subset of the facets that contain $f$.
    There the first entry is a single face, that is, one connected component of an intersection $F_{\cap J}$ of facets, while here it is the set $J$ itself.
    This loses nothing, since the refined subproduct indexed by $(J, L)$ is the disjoint union of its restrictions to the connected components of $F_{\cap J}$.
\end{remark}

For any $(J, L) \leq (J', L')$, by considering the inclusion $\zeta|_{F_{\cap J}} \subseteq \zeta|_{F_{\cap J'}}$ and the inclusion $(\mathrm{D}^2, \mathrm{S}^1)^L \subseteq (\mathrm{D}^2, \mathrm{S}^1)^{L'}$, we have a natural inclusion
\[
    ({\mathrm{D}^2}, \mathrm{S}^1)^{(\zeta, J, L)} \hookrightarrow ({\mathrm{D}^2}, \mathrm{S}^1)^{(\zeta, J', L')} .
\]

\begin{proposition}\label{prop:CWdecomp}
    The space $({\mathrm{D}^2}, \mathrm{S}^1)^{\zeta}$ admits a cell-decomposition in which all the refined subproducts
    \[
        ({\mathrm{D}^2}, \mathrm{S}^1)^{(\zeta, J, L)},
        \qquad L\subseteq J\subseteq[m],
    \]
    are CW-subcomplexes. 
\end{proposition}

\begin{proof}
    We may replace $\zeta$ by an isomorphic smooth principal bundle, so that $E:=\zeta\times({\mathrm{D}^2})^m$ is a smooth manifold with corners. The  facets of $E$ are
    \[
        H_j:=\zeta|_{F_{j}}\times({\mathrm{D}^2})^m , \qquad V_i:=\zeta\times \mathrm{S}^1_{(i)}\times\prod_{k\neq i}{\mathrm{D}^2_{(k)}}, \quad i, j\in[m].
    \]

    According to \cite{MR500991}, we can choose a triangulation of the manifold with corners $E$ which is compatible with its faces. That is, there exists a simplicial complex $K$ and a homeomorphism $|K|\to E$ such that the inverse image of each closed face is the realization of a subcomplex of $K$. For any $L\subseteq J\subseteq[m]$, the refined subproduct is the intersection of the corresponding facets,
    \[
        ({\mathrm{D}^2}, \mathrm{S}^1)^{(\zeta, J, L)}=\bigcap_{j\in J}H_j\cap\bigcap_{i\notin L}V_i ,
    \]
    and hence is a face of $E$, so it is a CW-subcomplex of the chosen triangulation. Their union is $({\mathrm{D}^2}, \mathrm{S}^1)^{\zeta}$, which is therefore also a CW-subcomplex, and the induced cell-decomposition has the required property.
\end{proof}

\begin{definition}\label{def:diagram}
    Let $Q$ and $\zeta\to Q$ be as in \Cref{def:genpolyprod}. 
    Let $\mathbf D:\mathscr P\to\mathbf{Top}$ be the diagram defined by
    \[
        \mathbf D(J, L):=({\mathrm D^2}, \mathrm S^1)^{(\zeta, J, L)}.
    \]
    For $(J, L)\leq(J', L')$, the structure map $d_{(J, L), (J', L')}:\mathbf D(J, L)\hookrightarrow\mathbf D(J', L')$ is the inclusion.
\end{definition}
Clearly, $(\mathrm D^2, \mathrm S^1)^{\zeta}$ is the colimit of the diagram $\mathbf D$. Each object of $\mathbf D$ carries the action of $\T^{n+m}$ given by \eqref{eq:actonassociatedbundle}, and the structure maps of $\mathbf D$ are equivariant with respect to these actions.

\begin{proposition}\label{prop:generalcollapsing}
    The collapsing map
    \[
        \hocolim_{\mathscr P} \mathbf D\stackrel{\simeq}{\longrightarrow} \colim_{\mathscr P} \mathbf D =  (\mathrm D^2, \mathrm S^1)^{\zeta}
    \]
    is an equivariant weak equivalence. We use $\simeq$ to denote weak equivalence.
\end{proposition}
\begin{proof}
    By \Cref{prop:CWdecomp}, each latching map of this diagram is the inclusion of a CW-subcomplex. 
    It follows from \Cref{lemma:latchinglemma} that the collapsing map is a weak equivalence. Moreover, the collapsing map is equivariant by its construction.
\end{proof}

\subsection{The Borel construction as a homotopy colimit}
We ultimately need to compute the equivariant cohomology, which is the cohomology of the Borel construction.
The following lemma lets us pass from weak equivalences of $G$-spaces to weak equivalences of their Borel constructions.
\begin{lemma}\label{lemma:equiinduceequiforbor}
    Let $G$ be a group.
    Assume $f:X\to Y$ is a $G$-equivariant weak equivalence between $G$-spaces. Then $f$ induces a weak equivalence from $X_G$ to $Y_G$.
\end{lemma}

\begin{proof}
    Consider the following diagram of fiber bundles:
    \[\begin{tikzcd}
	G & {EG\times X} & {X_G} \\
	G & {EG\times Y} & {Y_G}
	\arrow[from=1-1, to=1-2]
	\arrow["{=}", from=1-1, to=2-1]
	\arrow[from=1-2, to=1-3]
	\arrow["{\id\times f}", from=1-2, to=2-2]
	\arrow["{\id \times_{G} f}", from=1-3, to=2-3]
	\arrow[from=2-1, to=2-2]
	\arrow[from=2-2, to=2-3]
\end{tikzcd}\]
The conclusion follows by considering the long exact sequence of homotopy groups and then applying the five lemma.
\end{proof}
We write $\mathbf D_{\T^{n+m}}$ or $E\T^{n+m} \times_{\T^{n+m}} \mathbf{D}$ for the diagram obtained by taking Borel constructions objectwise. Since the structure maps of $\mathbf D$ are $\T^{n+m}$-equivariant, they induce maps between these Borel constructions, so $\mathbf D_{\T^{n+m}}$ is again a diagram over $\mathscr P$. 
One checks the natural identifications
\begin{equation}\label{eq:commutofcolimandBorel}
        (\colim_{\mathscr P}\mathbf D)_{\T^{n+m}}\cong\colim_{\mathscr P}(\mathbf D_{\T^{n+m}}),
    \qquad (\hocolim_{\mathscr P}\mathbf D)_{\T^{n+m}}\cong\hocolim_{\mathscr P}(\mathbf D_{\T^{n+m}}).
\end{equation}
Combining \Cref{{prop:generalcollapsing}} with \Cref{lemma:equiinduceequiforbor}, we have 
\begin{equation*}
    (\hocolim_{\mathscr P}\mathbf D)_{\T^{n+m}} \simeq (\colim_{\mathscr P}\mathbf D)_{\T^{n+m}}.
\end{equation*}
Together with the natural identifications in \eqref{eq:commutofcolimandBorel}, this gives the following lemma:
\begin{lemma}\label{lem:weakequivalenceOfhocolimandpolyhedralProduct}
    For the diagram $\mathbf{D}$, there is a weak equivalence:
    \begin{equation*}
        \hocolim_{\mathscr{P}} (\mathbf{D}_{\T^{n+m}} ) \simeq (\colim_{\mathscr P}\mathbf D)_{\T^{n+m}} = E\T^{n+m} \times_{\T^{n+m}} (\mathrm D^2, \mathrm S^1)^{\zeta}.
    \end{equation*}
    For each $(J, L) \in \mathscr{P}$, the natural projection 
    \begin{equation*}
        \mathbf{D}(J, L)_{\T^{n+m}} \to B\T^{n+m}
    \end{equation*}
    induces a commutative diagram:
    \[\begin{tikzcd}
        {\hocolim_{\mathscr{P}} (\mathbf{D}_{\T^{n+m}} )} && E\T^{n+m} \times_{\T^{n+m}} (\mathrm D^2, \mathrm S^1)^{\zeta} \\
        & {B\T^{n+m}}
        \arrow["\simeq", from=1-1, to=1-3]
        \arrow[from=1-1, to=2-2]
        \arrow[from=1-3, to=2-2]
    \end{tikzcd}\]
\end{lemma}

By \Cref{thm:mamehtpp}, $\mathscr{L}_{(Q, \zeta)}$ is $\T^{n+m}$-equivariantly homeomorphic to $(\mathrm{D}^2, \mathrm{S}^1)^{\zeta}$. 
According to \Cref{lem:weakequivalenceOfhocolimandpolyhedralProduct}, to compute the equivariant cohomology ring and $H^*(B\T^{n+m})$-module structure of $\mathscr{L}_{(Q, \zeta)}$, we only need to compute those of $\hocolim_{\mathscr{P}} (\mathbf{D}_{\T^{n+m}} )$.  

\subsection{The diagram of cohomology rings}\label{subsec:cohomologydiagram}
It remains to compute $H^*(\hocolim_{\mathscr{P}} (\mathbf{D}_{\T^{n+m}} ))$ together with its $H^*(B\T^{n+m})$-module structure. We therefore study the diagram $H^*(\mathbf{D}_{\T^{n+m}})$, which is obtained by applying the cohomology functor $H^*(-;\Z)$ to $\mathbf{D}_{\T^{n+m}}$ objectwise. Its structure maps are the ring homomorphisms between cohomology rings induced by the structure maps of $\mathbf{D}_{\T^{n+m}}$, and it is a diagram of $H^*(B\T^{n+m})$-modules. This diagram plays an important role in the computation of $H^*(\hocolim_{\mathscr{P}} (\mathbf{D}_{\T^{n+m}} ))$ and its $H^*(B\T^{n+m})$-module structure (see \Cref{lemma:spectralsequence}). 

For each $(J, L)\in \mathscr{P}$, we have 
\begin{equation}\label{eq:defofDT}
    \begin{aligned}
        \mathbf{D}_{\T^{n+m}} (J, L) &= E\T^{n+m} \times_{\T^{n+m}} (\mathrm{D}^2, \mathrm{S}^1)^{(\zeta, J, L)} \\
        &= E\T^{n+m} \times_{\T^{n+m}} \bigl(\zeta|_{F_{\cap J}}\times\prod_{j\in L}{\mathrm{D}^2_{(j)}}\times\prod_{j\notin L}\mathrm{S}^1_{(j)} , \bigr) \\ 
        &= (E\T^n \times_{\T^n} \zeta|_{F_{\cap J}})\times\prod_{j\in L}(E\T^1 \times_{\T^1}{\mathrm{D}^2_{(j)}})\times\prod_{j\notin L}(E\T^1 \times_{\T^1} \mathrm{S}^1_{(j)}).
    \end{aligned}
\end{equation}
The projection map $ \mathbf{D}_{\T^{n+m}} (J, L) \to B\T^{n+m}$ is then induced by the projection map for each term among the above product:
\begin{equation*}
    E\T^n \times_{\T^n} \zeta|_{F_{\cap J}} \to B\T^n, \quad E\T^1 \times_{\T^1}{\mathrm{D}^2_{(j)}} \to B\T^1, \quad E\T^1 \times_{\T^1} \mathrm{S}^1_{(j)} \to B\T^1. 
\end{equation*}
Also, for $(J, L) \leq (J', L')$, the structure map $\mathbf{D}(J, L) \to  \mathbf{D}(J', L')$ is induced by 
\begin{equation*}
    \begin{aligned}
        E\T^n \times_{\T^n} \zeta|_{F_{\cap J}} &\to E\T^n \times_{\T^n} \zeta|_{F_{\cap J'}} \\
    E\T^1 \times_{\T^1} {\mathrm{S}^1_{(j)}} &\to E\T^1 \times_{\T^1} {\mathrm{D}^2_{(j)}}, \quad j\in L' \setminus L, 
    \end{aligned}
\end{equation*}
and be the identity map for the remaining term of the product. 

Since $\T^{n+m}=\T^n\times\T^m$ acts on the two factors of \eqref{eq:defofDT} separately, each object of $\mathbf D_{\T^{n+m}}$ splits as a product of
\[
    (\zeta|_{F_{\cap J}})_{\T^n}=E\T^n\times_{\T^n}\zeta|_{F_{\cap J}}
    \qquad\text{and}\qquad
    (\mathrm{D}^2,\mathrm{S}^1)^{L}_{\T^m}=E\T^m\times_{\T^m}(\mathrm{D}^2,\mathrm{S}^1)^{L},
\]
the first a space over $B\T^n$ and the second a space over $B\T^m$, and the structure maps split accordingly. The cohomology ring of $(\mathrm{D}^2,\mathrm{S}^1)^{L}_{\T^m}$ together with its $H^*(B\T^m)$-module structure is standard, see \cite[Section~4.3]{MR3363157}. To compute those of $(\zeta|_{F_{\cap J}})_{\T^n}$, we use the following fact on the Borel construction of a principal bundle, which we take from \cite[p.~60]{MR889050}. 
\begin{lemma}\label{lem:tomDieck}
    Let $P\to B$ be a principal $G$-bundle, the projection $r: EG\times_G P \to B$ is a homotopy equivalence, and there is a pullback diagram:
\[\begin{tikzcd}
	P & EG \\
	B & BG
	\arrow["{f_*}", from=1-1, to=1-2]
	\arrow["{\pi_P}"', from=1-1, to=2-1]
	\arrow["{\pi_{EG}}", from=1-2, to=2-2]
	\arrow["f", from=2-1, to=2-2]
\end{tikzcd}, \]
where $f_*$ is $G$-equivariant.
With this diagram, a homotopy inverse of $r$ is given by $s: B \to EG \times_G P$:
\begin{equation*}
    s(x) = [ f_*(v), v], \quad \forall \ v\in \pi_P^{-1}(x).
\end{equation*}
The composition $s$ with the projection $EG\times_G P\to BG$ equals $f$.     
\end{lemma}

\begin{lemma}\label{lem:cohomologyofBorel}
    For each $J\subseteq[m]$, there is an isomorphism
    \[
        H^*\bigl((\zeta|_{F_{\cap J}})_{\T^n}\bigr)\cong H^*(F_{\cap J}).
    \]
    Under this isomorphism, the generator $y_i$ of $H^*(B\T^n)=\Z[y_1,\ldots,y_n]$ acts by cup product with $c(\zeta)_i|_{F_{\cap J}}$, that is,
    \[
        y_i\cdot\alpha=c(\zeta)_i|_{F_{\cap J}}\cup\alpha,
        \qquad \alpha\in H^*(F_{\cap J}).
    \]
    Furthermore, for $J'\subseteq J$, the homomorphism
    \[
        H^*\bigl((\zeta|_{F_{\cap J'}})_{\T^n}\bigr)\longrightarrow H^*\bigl((\zeta|_{F_{\cap J}})_{\T^n}\bigr)
    \]
    under these isomorphisms, the restriction map $H^*(F_{\cap J'})\to H^*(F_{\cap J})$ induced by the inclusion $F_{\cap J}\hookrightarrow F_{\cap J'}$.
\end{lemma}
\begin{proof}
    For $J' \subseteq J \subseteq [m]$, here is a commutative diagram by \Cref{lem:tomDieck}:
\[\begin{tikzcd}
	{F_{\cap J}} & {(\zeta|_{F_{\cap J}})_{\T^n}} & {B\T^n} \\
	{F_{\cap J'}} & {(\zeta|_{F_{\cap J'}})_{\T^n}} & {B\T^n}
	\arrow["\simeq", from=1-1, to=1-2]
	\arrow[hook, from=1-1, to=2-1]
	\arrow[from=1-2, to=1-3]
	\arrow[from=1-2, to=2-2]
	\arrow["{=}", from=1-3, to=2-3]
	\arrow["\simeq", from=2-1, to=2-2]
	\arrow[from=2-2, to=2-3]
\end{tikzcd}\]
where the two horizontal arrows in the left square are the sections given by \Cref{lem:tomDieck}.  
Let $\chi: Q\to B\T^n$ be the classifying map of $\zeta$. The two horizontal composites are
\begin{equation}\label{eq:compositionisclassifyingmap}
    \chi|_{F_{\cap J}}\colon F_{\cap J}\to B\T^n
    \qquad\text{and}\qquad
    \chi|_{F_{\cap J'}}\colon F_{\cap J'}\to B\T^n.
\end{equation} 
Applying the cohomology functor, the assertions of the lemma follow directly from this diagram together with \eqref{eq:compositionisclassifyingmap}. 
\end{proof}

Now we can compute the ring and $H^*(B\T^{n+m})$-module structure of the diagram $H^*(\mathbf{D}_{\T^{n+m}})$. 
For the diagram $\mathbf{D}_{\T^{n+m}}$, we have 
\begin{equation*}
    \begin{aligned}
        H^*(\mathbf{D}_{\T^{n+m}})(J, L) &:= H^*(\mathbf{D}_{\T^{n+m}}(J, L)) \\ 
        &=H^*_{\T^n}( \zeta|_{F_{\cap J}}) \otimes \bigotimes_{j\in L}H^*_{\T^1} (\mathrm{D}^2_{(j)}) \otimes 
        \bigotimes_{j \notin L} H^*_{\T^1} (\mathrm{S}^1_{(j)}) \\ 
        &=H^*(F_{\cap J}) \otimes \Z[ \mathbf{x}_{L}]. 
    \end{aligned}
\end{equation*}
For $(J, L) \leq  (J', L')$, 
the structure map $\varepsilon_{(J', L'), (J, L)}: H^*(\mathbf{D}_{\T^{n+m}})(J', L') \to H^*(\mathbf{D}_{\T^{n+m}})(J, L)$ is given by 
        \begin{equation}\label{eq:varepsilonG}
                \alpha\otimes f((x_i)_{i\in L'}) \mapsto(\alpha|_{F_{\cap J}}\otimes1)f((u_i)_{i\in L'}),
        \end{equation}
        for $\alpha \in H^*(F_{\cap J'})$ and $f \in \Z[\mathbf{x}_{L'}]$, 
        where 
        \begin{equation*}
            u_i=
                    \begin{cases}
                        1\otimes x_i  \quad \text{if } i\in L, \\ 
                        0 \quad \text{if } i\in L'\setminus L.
                    \end{cases}
        \end{equation*} 
The $H^*(B\T^{n+m})$-module structure, under the isomorphism $H^*(B\T^{n+m}) = \Z[\mathbf{y}_{[n]}, \mathbf{x}_{[m]}]$, is given by
\begin{equation}\label{def:xiactiondefinemodule}
    \begin{aligned}
        x_i\cdot(\alpha\otimes f)&=
        \begin{cases}
            0 &\text{if }i\notin L, \\
            \alpha\otimes x_i f &\text{if }i\in L,
        \end{cases}\\
        y_k\cdot(\alpha\otimes f)&=(c(\zeta)_k|_{F_{\cap J}}\cup\alpha)\otimes f.
    \end{aligned}
\end{equation}

    \section{Equivariant cohomology rings and module structures over $H^*(B\T^{n+m})$}\label{sec:computation}
In this section we compute $H^*(\hocolim_{\mathscr P}\mathbf D_{\T^{n+m}};\Z)$ together with its $H^*(B\T^{n+m})$-module structure, and deduce \Cref{thm:generalpolyhedral} and \Cref{THM:1}.

\subsection{The limit and the topological face ring}
To simplify notation, we write
\[
    \mathbf G:=H^*(\mathbf D_{\T^{n+m}};\Z)
\]
for the diagram of cohomology rings computed in \Cref{subsec:cohomologydiagram}. Thus $\mathbf G$ is a diagram over $\mathscr P^{op}$ with
\[
    \mathbf G(J,L)=H^*(F_{\cap J};\Z)\otimes\Z[\mathbf x_L],
\]
whose structure maps are the homomorphisms $\varepsilon_{(J',L'),(J,L)}$ of \eqref{eq:varepsilonG}, and which is a diagram of $\Z[\mathbf x_{[m]},\mathbf y_{[n]}]$-modules with the action \eqref{def:xiactiondefinemodule}. Since its objects are cohomology rings of spaces, $\mathbf G$ is a diagram of graded commutative rings. For $q\ge0$, we define a diagram $\mathbf G^q$ over $\mathscr P^{op}$ as follows. For each $(J,L)\in\mathscr P$, let $\mathbf G^q(J,L)$ be the subgroup of $\mathbf G(J,L)$ consisting of the elements of degree $q$. Since the structure maps of $\mathbf G$ preserve the degree, their restrictions give the structure maps of $\mathbf G^q$. Then
\[
    \mathbf G=\bigoplus_{q\ge0}\mathbf G^q.
\]

The limit $\lim_{\mathscr{P}^{op}}\mathbf{G}$ is, by definition, the subring of $\bigoplus_{(J, L) \in \mathscr{P}} \mathbf{G}(J, L)$ consisting of those families $b=(b_{(J, L)})_{(J, L) \in \mathscr{P}}$ which satisfy
\begin{equation}\label{eq:limitcondition}
    \varepsilon_{(J', L'), (J, L)}(b_{J', L'}) = b_{(J, L)} \qquad \text{whenever } (J, L)\le (J', L').
\end{equation}
We call $b_{(J, L)}$ the $(J, L)$-component of the element $b$ in the limit. 

We shall identify the limit with the topological face ring $\Z\langle Q\rangle$ of \Cref{def:topologicalfacering}, whose $\star$-multiplication is defined there, and we endow it with the $c(\zeta)$-twisted module structure of \Cref{def:twistedZxmmod}.

\begin{definition}\label{def:twistedZxmmod}
    For $d\in H^2(Q;\Z^n)$, we call the module structure on $\Z\langle Q\rangle$ defined by \eqref{eq:intro-module-one} the \emph{$d$-twisted} $\Z[\mathbf x_{[m]}, \mathbf y_{[n]}]$-module structure.
\end{definition}

\begin{lemma}\label{prop:limascompatible}
    The limit $\lim_{\mathscr P^{op}}\mathbf G$ is canonically isomorphic to the ring of families
    \[
        \mathcal H=\left\{(b_J)_{J\subseteq[m]}\in\bigoplus_{J\subseteq[m]}\mathbf G(J, J)\ \middle|\
        \begin{array}{l}
            \varepsilon_{(J, J), (J, L)}(b_J)=\varepsilon_{(L, L), (J, L)}(b_L)\\
            \text{for every }L\subseteq J\subseteq[m]
        \end{array}\right\}.
    \]
    We call $b_J$ the $J$-component of $b\in\mathcal H$.
\end{lemma}
\begin{proof}
    For each $(J, L)\in\mathscr P$, the diagonal object $(L, L)$ belongs to $\mathscr P$ and lies above $(J, L)$. By the compatibility condition given in \eqref{eq:limitcondition}, an element of the limit is determined by its diagonal components which satisfy the stated condition.
    Conversely, a family in $\mathcal H$ determines all components by
    \[
        b_{(J, L)}:=\varepsilon_{(L, L), (J, L)}(b_L).
    \]
    For an arrow $(J, L)\leq(J', L')$, use the compatibility for $L\subseteq L'$ at $(L', L)$, then restrict to $F_{\cap J}$. This gives the limit condition. These constructions are mutually inverse.
\end{proof}

\begin{proposition}\label{prop:gradedringiso}
    There is an isomorphism of graded rings:
    \begin{equation}
        \lim\nolimits_{\mathscr{P}^{op}} \mathbf{G} \cong \Z\langle Q\rangle.
    \end{equation}
\end{proposition}

To prove this proposition, we fix the following notation for \emph{monomials}.
A monomial is an expression
\[
    \mu=\prod_{i\in[m]}x_i^{a_i}, \qquad a_i\ge0.
\]
We write
\[
    |\mu|=\sum_{i\in[m]}a_i, \qquad \operatorname{supp}\mu=\{i\in[m]:a_i>0\}.
\]
Thus $\mu$ has degree $2|\mu|$.
The decomposition of an element of $H^*(F_{\cap J};\Z)\otimes\Z[\mathbf x_L]$ as a sum of tensor products of classes in $H^*(F_{\cap J};\Z)$ with distinct monomials plays an important role in the following proof and in the next subsection.

\begin{proof}
    According to the preceding lemma, to prove this proposition, we only need to construct a ring homomorphism $\rho: \Z\langle Q\rangle \to \mathcal{H}$ and show it is a ring isomorphism. Explicitly, we define $\rho$ on a term $\alpha \otimes f \in H^*(F_{\cap J};\Z) \otimes \Z[\mathbf{x}^+_J]$ by
    \begin{equation}
        (\rho(\alpha \otimes f))_{J'} = 
        \begin{cases}
            \alpha|_{F_{\cap J'}} \otimes f  \quad \text{if }J\subseteq J'; \\ 
            0 \quad \quad  \quad \quad \quad \quad \text{otherwise}.
        \end{cases}
    \end{equation}
    One can check that $\rho: \Z\langle Q \rangle \to \mathcal{H}$ is a ring homomorphism.

    Injectivity: Suppose that $\rho(a)=0$ for some $a=(a_J) \in \bigoplus_{J\subseteq[m]}H^*(F_{\cap J}) \otimes \Z[\mathbf{x}^+_J]$ in $\Z\langle Q\rangle$.
    If $a \neq 0$, we choose $J\subseteq[m]$ which is minimal with respect to inclusion among those with $a_J\neq 0$. Thus $a_L=0$ for every $L \subsetneq J$. By the definition of $\rho$, we have $(\rho(a))_J=a_J=0$. This is a contradiction. Thus we have $a=0$. 

    Surjectivity: Write each component uniquely as a finite sum
    $b_J=\sum_{\operatorname{supp}\mu\subseteq J}\beta_{J, \mu}\otimes\mu$
    over monomials $\mu$, where $\beta_{J, \mu}\in H^*(F_{\cap J};\Z)$.
    For $S=\operatorname{supp}\mu\subseteq J$, the compatible condition for elements in $\mathcal{H}$ gives
    \begin{equation*}
        \varepsilon_{(J, J), (J, S)}(b_J)=\varepsilon_{(S, S), (J, S)}(b_S)
    \end{equation*}
    Since the structure map $\varepsilon_{(J, J), (J, S)}$ sets the variables outside $S$ to zero, this gives
    \begin{equation}\label{eq:coherent}
        \beta_{J, \mu}=\beta_{S, \mu}|_{F_{\cap J}}.
    \end{equation}
    Hence $a:=\sum_\mu\beta_{\operatorname{supp}\mu, \mu}\otimes\mu$ belongs to $\Z\langle Q\rangle$ and satisfies $\rho(a)=b$. 
\end{proof}

\subsection{Higher limits and the cohomology ring}\label{subsec:higherlimits}
In \cite{MR4721373}, Yu treats the case $n=0$ and computes the cohomology and equivariant cohomology of $(\mathrm D^2, \mathrm S^1)^{\zeta}$ by a stable decomposition.
In that case each refined subproduct $({\mathrm D^2}, \mathrm S^1)^{(\zeta, J, L)}$ is the product $F_{\cap J}\times\prod_{j\in L}\mathrm D^2_{(j)}\times\prod_{j\notin L}\mathrm S^1_{(j)}$, so after one suspension it splits as a wedge of smash products (see \cite[Section~2.4]{MR4721373}).
We compute the cohomology of $\hocolim_{\mathscr P}\mathbf{D}_{\T^{n+m}}$ instead by tools from homotopy colimits.

We have identified $\Z\langle Q\rangle$ with $\lim_{\mathscr P^{op}}\mathbf G$ by \Cref{prop:gradedringiso}.
To relate this limit to the cohomology of $\hocolim_{\mathscr P}\mathbf{D}_{\T^{n+m}}$, we use the following spectral sequence.
\begin{lemma}[{{\cite[Theorem 18.3]{dugger2008primer}, \cite[Proposition 5.1]{MR232393}}}]\label{lemma:spectralsequence}
    Let $D: I\to \mathbf{Top}$ be a diagram. 
    There is a spectral sequence $E_2^{p, q}=\lim\nolimits^p_{I^{op}} H^q(D; \Z)\Rightarrow  H^{p+q}(\mathrm{hocolim}\, D, \Z)$. The differentials have the form $d_r\colon E_r^{p, q}\to E_r^{p+r, q-r+1}$.
\end{lemma}

Here diagram cohomology is expressed in terms of derived limits as in \cite[Section~17]{dugger2008primer}. We shall use the following facts about these functors.
\begin{enumerate}
    \item The \(0\)-th derived limit is the usual limit functor. In our case,
    \[
        \lim\nolimits^0_{\mathscr{P}^{op}} \mathbf{G}
        = \lim\nolimits_{\mathscr{P}^{op}} \mathbf{G}.
    \]
    \item One way to define the derived limit
    \(\lim\nolimits^n_{\mathscr{P}^{op}} \mathbf{G}\)
    is as the cohomology group of a cochain complex whose \(n\)-th cochain group
    \(C^n(\mathscr{P}^{op}, \mathbf{G})\) is given by
    \begin{equation}
        C^n(\mathscr{P}^{op}, \mathbf{G})
        = \bigoplus_{i_0 \to i_1 \to \cdots \to i_n} \mathbf{G}(i_n)
    \end{equation}
    (see {\cite[Definition 3.1]{husainov2023criteriapreservingcategorycohomology}}),
    where \(i_0, \dots, i_n\) are objects in \(\mathscr{P}^{op}\), and each arrow represents a morphism in \(\mathscr{P}^{op}\). Note that the partial order relation among \(i_0, \dots, i_n\) in \(\mathscr{P}\) is
    \[
        i_n \leq i_{n-1} \leq \cdots \leq i_0.
    \]
    \item Direct sums of diagrams are taken objectwise. The cochain description shows that finite direct sums commute with derived limits.
\end{enumerate}

For the reader who wants to know more about the derived limit functors $\lim\nolimits^n_{I}$, we refer to \cite[Section 17]{dugger2008primer} and \cite[Section 3.5]{MR1269324} for more details.

Our goal is to show that $\lim\nolimits^n_{\mathscr{P}^{op}} \mathbf{G}=0$ for $n>0$ (see \Cref{cor:vanishingG}).
To prove this vanishing result, we introduce some auxiliary functors. 
Fix a covariant functor $G\colon 2^{[m]}\to\mathbf{Ab}$, where $2^{[m]}$ denotes the power set of $[m]$ ordered by inclusion. 
For a fixed integer $d\geq 0$, define a functor $A(G, d)\colon\mathscr{P}^{op}\to\mathbf{Ab}$ by
\[
    A(G, d)(J, L) \;=\; G(J)\otimes\Z[\mathbf{x}_L]_{2d},
\]
where $\Z[\mathbf{x}_L]_{2d}$ denotes the homogeneous degree-$2d$ part of the polynomial ring $\Z[\mathbf{x}_L]$ (note that $\deg x_i=2$ for each variable). 
For $(J, L)\le(J', L')$, the structure map $A(G, d)(J', L')\to A(G, d)(J, L)$ is the tensor product of the homomorphism $G(J')\to G(J)$ with
\[
    \rho_{L, L'}\colon\Z[\mathbf{x}_{L'}]_{2d}\to\Z[\mathbf{x}_L]_{2d}, \qquad
    x_i\mapsto
    \begin{cases}
        x_i, & i\in L, \\
        0, & i\in L'\setminus L.
    \end{cases}
\]

The following lemma gives a decomposition of $A(G, d)$ according to monomials of $\Z[\mathbf{x}_L]_{2d}$. 
\begin{lemma}\label{lem:monomialbasis}
    For each $L\subseteq[m]$ there is a direct sum decomposition
    \[
        \Z[\mathbf{x}_L]_{2d} \;=\; \bigoplus_{\substack{\mu:\, |\mu|=d\\ \operatorname{supp}\mu\subseteq L}} \Z[\mu],
    \]
    compatible with the maps $\rho_{L, L'}$.
    Here, we use $\Z[\mu]$ to denote the subgroup of $\Z[\mathbf{x}_L]$ generated by $\mu$. 
    Consequently, there is a decomposition of the diagram $A(G, d)$:
    \[
        A(G, d) = \bigoplus_{|\mu|=d} A_{\mu}, \quad
        A_{\mu}(J, L) := G(J)\otimes\mathbf{1}_{\mu}(L),
    \]
    where $\mathbf{1}_{\mu}(L)=\Z$ if $\operatorname{supp}\mu\subseteq L$ and $\mathbf{1}_{\mu}(L)=0$ otherwise.
\end{lemma}
\begin{proof}
    The degree-$2d$ part of $\Z[\mathbf{x}_L]$ is spanned by the monomials $\mu$ with $|\mu|=d$. Under $\rho_{L, L'}$ such a $\mu$ is sent to itself if $\operatorname{supp}\mu\subseteq L$ and to $0$ otherwise.
\end{proof}

The following lemmas will allow us to restrict to these objects and then simplify the indexing category.
\begin{lemma}[Extension by zero]\label{lem:extensionbyzero}
    Let $U\subseteq\mathscr{P}$ be an \emph{up-set}, i.e. $u\in U$ and $u\le v$ imply $v\in U$, and let $D\colon U^{op}\to\mathbf{Ab}$ be a functor. 
    The extension by zero $eD\colon\mathscr{P}^{op}\to\mathbf{Ab}$, given by $(eD)(u)=D(u)$ for $u\in U$ and $(eD)(u)=0$ otherwise, satisfies
    \[
        \lim\nolimits^{n}_{\mathscr{P}^{op}}(eD) \;\cong\; \lim\nolimits^{n}_{U^{op}}D \qquad (n\ge0).
    \]
\end{lemma}
\begin{proof}
    The cochain group $C^n(\mathscr{P}^{op};eD)$ is given by
    $\bigoplus_{u_0\le\cdots\le u_n}(eD)(u_0)$, the indexed chain $u_0\le\cdots\le u_n$ being taken in $\mathscr P$.
    Since $eD$ vanishes outside $U$, a summand $(eD)(u_0)$ is nonvanishing only if $u_0\in U$.
    Because $U$ is an up-set and $u_0\le\cdots\le u_n$, this implies $u_1, \cdots, u_n\in U$ as well.
    Thus, by dropping these vanishing summands $(eD)(u_0)$ with $u_0 \notin U$, we obtain summands indexed by chains entirely contained in $U$.
    Therefore, we can identify the cochain complex $C^n(\mathscr{P}^{op};eD)$, including its differential, with $C^n(U^{op};D)$. This completes the proof.
\end{proof}

The comma category is a basic notion of category theory, which we shall use below. We recall its definition.

\begin{definition}\label{def:commacategory}
    The \emph{comma category} $(F\downarrow d)$ of a functor $F\colon I_1\to I_2$ over an object $d\in I_2$ has as objects the pairs $(c, f)$, where $c\in I_1$ and $f\colon F(c)\to d$ is a morphism of $I_2$, and as morphisms $(c, f)\to(c', f')$ the morphisms $g\colon c\to c'$ of $I_1$ with $f=f'\circ F(g)$ (see the following diagram). 
    \[\begin{tikzcd}
	{F(c')} && {F(c)} \\
	& d
	\arrow["{f'}"', from=1-1, to=2-2]
	\arrow["{F(g)}"', from=1-3, to=1-1]
	\arrow["f", from=1-3, to=2-2]
\end{tikzcd}\]
    If there is no morphism in $I_2$ from $F(c)$ to $d$ for any $c\in I_1$, the comma category $(F\downarrow d)$ is empty.
\end{definition}

All the categories that appear in this section are induced by posets. In that case the comma category has the following simple description.

\begin{lemma}\label{lem:commaposet}
    Let $I_1$ and $I_2$ be categories induced by posets and let $F\colon I_1\to I_2$ be a functor. For every object $d$ of $I_2$, the comma category $(F\downarrow d)$ is the full subcategory of $I_1$ on the objects $c$ with $F(c)\le d$.
\end{lemma}

\begin{proof}
    In a category induced by a poset there is at most one morphism between two objects. By \Cref{def:commacategory}, an object of $(F\downarrow d)$ is a pair $(c, f)$ with $f\colon F(c)\to d$ a morphism of $I_2$. Such a morphism exists only when $F(c)\le d$, and then it is unique, so we may drop $f$ and regard the objects of $(F\downarrow d)$ as the objects $c$ of $I_1$ with $F(c)\le d$. A morphism $(c, f)\to(c', f')$ of $(F\downarrow d)$ is a morphism $g\colon c\to c'$ of $I_1$ with $f=f'\circ F(g)$. By the same uniqueness, such a $g$ exists only when $c\le c'$, and in that case the identity $f=f'\circ F(g)$ holds automatically. Hence $(F\downarrow d)$ and the full subcategory of $I_1$ on these objects have the same objects and the same morphisms.
\end{proof}

\begin{remark}\label{rem:commaposetopp}
    The conclusion of \Cref{lem:commaposet} still holds when the source category is replaced by its opposite. Indeed, if $I_1$ is induced by a poset, then $I_1^{\mathrm{op}}$ is induced by the poset obtained by reversing the order relation. Hence a functor $F\colon I_1^{\mathrm{op}}\to I_2$ is a functor between categories induced by posets, and \Cref{lem:commaposet} applies with $I_1$ replaced by $I_1^{\mathrm{op}}$. The case of a functor $F\colon I_1\to I_2^{\mathrm{op}}$ is analogous, and we leave the statement to the reader.
\end{remark}

\begin{lemma}[Verdier criterion {{\cite[Corollary 4.6]{husainov2023criteriapreservingcategorycohomology}}}]\label{lem:finalfunctor}
    Let $F\colon I_1\to I_2$ be a functor such that, for every $d\in I_2$, the classifying space of the comma category  $(F\downarrow d)$ is non-empty and contractible. 
    Then for every functor $D\colon I_2\to\mathbf{Ab}$,
    \[
        \lim\nolimits^{n}_{I_1}(D\circ F) \;\cong\; \lim\nolimits^{n}_{I_2}D \qquad (n\ge0).
    \]
\end{lemma}

\begin{lemma}\label{lem:initialobject}
    If a category $I$ has an initial object, then $\lim\nolimits^{n}_{I} D=0$ for all $n\ge1$ and all functors $D: I\to\mathbf{Ab}$.
\end{lemma}

\begin{proof}
Let $i_0$ denote the initial object of $I$.  Let $I_0$ be the full subcategory of $I$ on the single object $i_0$, and let $F\colon I_0\to I$ be the inclusion.
Since $i_0$ is initial in $I$, for every $i\in I$ the comma category $(F\downarrow i)$ has exactly one object and exactly one morphism.
Thus \Cref{lem:finalfunctor} gives $\lim\nolimits^{n}_{I}D \cong \lim\nolimits^{n}_{I_0}(D\circ F)$.  For the second
term we have $\lim\nolimits^{0}_{I_0}(D\circ F)=D(i_0)$ and $\lim\nolimits^{n}_{I_0}(D\circ F)=0$ for $n\geq1$.
\end{proof}

For each monomial summand, the preceding lemmas reduce the computation to a category with an initial object.
This gives the following vanishing theorem.
\begin{theorem}\label{thm:vanishinghigherlimits}
    For every covariant functor $G\colon 2^{[m]}\to\mathbf{Ab}$ and every $d\ge0$,
    \[
        \lim\nolimits^{k}_{\mathscr{P}^{op}}A(G, d) \;=\; 0 \qquad (k\ge1).
    \]
\end{theorem}
\begin{proof}
    By \Cref{lem:monomialbasis} we have $A(G, d)=\bigoplus_{|\mu|=d}A_\mu$, a finite direct sum. Since $\lim\nolimits^{k}$ commutes with finite direct sums, it suffices to prove $\lim\nolimits^{k}_{\mathscr{P}^{op}}A_\mu=0$ for each monomial $\mu$ with $|\mu|=d$.
    Fix such a $\mu$ and put
    \[
        U_\mu \;=\; \{(J, L)\in\mathscr{P}:\operatorname{supp}\mu\subseteq L\}, \qquad
        V_{\mu} \;=\; \{J\subseteq[m]:\operatorname{supp}\mu\subseteq J\}.
    \]
    Let $Pr\colon U_\mu^{op}\to V_{\mu}$ be the projection $(J, L)\mapsto J$. 
    Note that $U_\mu$ is an up-set of $\mathscr{P}$ and that $A_\mu(J, L)=G(J)$ for $(J, L)\in U_\mu$ and $A_\mu(J, L)=0$ otherwise, i.e. $A_\mu$ is the extension by zero from $U_\mu$ of $G \circ Pr$. By \Cref{lem:extensionbyzero},
    \[
        \lim\nolimits^{k}_{\mathscr{P}^{op}}A_\mu \;\cong\; \lim\nolimits^{k}_{U_\mu^{op}}(G \circ Pr).
    \]

    By \Cref{lem:commaposet} and \Cref{rem:commaposetopp},
    for $J\in V_\mu$, the comma category $(Pr\downarrow J)$ can be identified with the full subcategory of $\mathscr P^{op}$: 
    \[
        \{(J', L'):\operatorname{supp}\mu\subseteq L'\subseteq J'\subseteq J\}.
    \] 
    For every object $(J', L')$, the relation
    \[
        (J, \operatorname{supp}\mu)\leq(J', L')\quad\text{in }\mathscr P
    \]
    gives a unique morphism $(J', L')\to(J, \operatorname{supp}\mu)$ in $\mathscr P^{op}$. Thus $(J, \operatorname{supp}\mu)$ is a terminal object of $(Pr\downarrow J)$, so its nerve complex is nonempty and contractible. By \Cref{lem:finalfunctor},
    \[
        \lim\nolimits^{k}_{U_\mu^{op}}(G \circ Pr) \;\cong\; \lim\nolimits^{k}_{V_{\mu}}G.
    \]
    Finally, $V_{\mu}$ has the initial object $\operatorname{supp}\mu$, so \Cref{lem:initialobject} gives $\lim\nolimits^{k}_{V_{\mu}}G=0$ for $k\ge1$.
\end{proof}

\begin{remark}\label{rem:commacube}
    The comma category also has a description in terms of cube faces. For $\operatorname{supp}\mu\subseteq L'\subseteq J'\subseteq J$, put
    \[
        C(J', L'):=\{t\in[0, 1]^m:t_i=0\text{ for }i\notin J', \ t_i=1\text{ for }i\in L'\}.
    \]
    Then $(Pr\downarrow J)$ is the poset of the nonempty faces $C(J', L')$ of $C(J, \operatorname{supp}\mu)$, ordered by inclusion. This poset contains the whole cube $C(J, \operatorname{supp}\mu)$ as its largest element. Hence its nerve complex is a cone with vertex $(J, \operatorname{supp}\mu)$, and it is the barycentric subdivision of $C(J, \operatorname{supp}\mu)$.
\end{remark}

\begin{corollary}\label{cor:vanishingG}
    For every $q\ge0$,
    \[
        \lim\nolimits^{k}_{\mathscr{P}^{op}} H^q(\mathbf{D}_{\T^{n+m}};\Z) \;=\; 0 \qquad (k\ge1).
    \]
\end{corollary}
\begin{proof}
    Since $\mathbf G^q$ is the degree-$q$ part of $\mathbf G=H^*(\mathbf D_{\T^{n+m}};\Z)$, it suffices to show that $\lim\nolimits^{k}_{\mathscr{P}^{op}}\mathbf{G}^q=0$ for $k\ge1$. As each $x_i$ has degree $2$,
    \[
        \mathbf{G}^q(J, L) \;=\; \bigoplus_{d=0}^{\lfloor q/2\rfloor} H^{q-2d}(F_{\cap J};\Z)\otimes\Z[\mathbf{x}_L]_{2d},
    \]
    and this decomposition is compatible with the structure maps. Writing $G_d$ for the covariant functor $J\mapsto H^{q-2d}(F_{\cap J};\Z)$ on $2^{[m]}$, we obtain $\mathbf{G}^q\cong\bigoplus_{d=0}^{\lfloor q/2\rfloor}A(G_d, d)$. Since $\lim\nolimits^{k}$ commutes with finite direct sums, \Cref{thm:vanishinghigherlimits} yields $\lim\nolimits^{k}_{\mathscr{P}^{op}}\mathbf{G}^q=0$ for $k\ge1$, which proves the claim.
\end{proof}
The vanishing of higher limits allows us to identify the cohomology of $\hocolim_{\mathscr P}\mathbf{D}_{\T^{n+m}}$ with the ordinary limit computed in the previous subsection.
\begin{corollary}\label{cor:spectralcollapse}
    The edge homomorphism induces graded ring isomorphisms which are $H^*(B\T^{n+m};\Z)$-linear
    \begin{equation}\label{eq:homcolimislim}
        H^*(\hocolim_{\mathscr{P}}\mathbf{D}_{\T^{n+m}};\Z)
        \cong\lim\nolimits_{\mathscr{P}^{op}}H^*(\mathbf{D}_{\T^{n+m}};\Z)
        \cong \Z\langle Q\rangle,
    \end{equation}
    where $\Z\langle Q\rangle$ carries the $c(\zeta)$-twisted module structure of \Cref{def:twistedZxmmod}.
\end{corollary}

\begin{proof}
    By \Cref{cor:vanishingG}, the spectral sequence of \Cref{lemma:spectralsequence} is concentrated in the column $p=0$ from the $E_2$-page onwards. Hence the edge homomorphism is an isomorphism. It is induced by restriction to the objects $\mathbf{D}_{\T^{n+m}}(J, L)$ over $B\T^{n+m}$, so it preserves cup products and is $H^*(B\T^{n+m};\Z)$-linear. The second isomorphism follows from the description of $\mathbf G$ in \Cref{subsec:cohomologydiagram} and \Cref{prop:gradedringiso}.
\end{proof}
We conclude this section with the following theorem, from which \Cref{THM:1} follows directly.
\begin{theorem}\label{thm:generalpolyhedral}
    There are graded ring isomorphisms
    \[
        H^*\bigl(((\mathrm D^2, \mathrm S^1)^{\zeta})_{\T^{n+m}};\Z\bigr)
        \cong \lim\nolimits_{\mathscr P^{op}}\mathbf G
        \cong \Z\langle Q\rangle.
    \]
    These are also isomorphisms of $H^*(B\T^{n+m};\Z)$-modules.
    On the limit, the module structure is induced by
    \eqref{def:xiactiondefinemodule} with the characteristic class
    $c(\zeta)$.
    On $\Z\langle Q\rangle$, it is given by $\star$-multiplication
    through the graded ring homomorphism
    \[
        x_i\longmapsto 1\otimes x_i,
        \qquad
        y_k\longmapsto c(\zeta)_k\otimes1.
    \]
\end{theorem}

\begin{proof}
    By \Cref{lem:weakequivalenceOfhocolimandpolyhedralProduct}, the Borel construction
    $((\mathrm D^2, \mathrm S^1)^{\zeta})_{\T^{n+m}}$
    is weakly equivalent to $\hocolim_{\mathscr P}\mathbf{D}_{\T^{n+m}}$
    over $B\T^{n+m}$.
    The result follows from \Cref{cor:spectralcollapse}.
\end{proof}

\begin{proof}[Proof of \Cref{THM:1}]
    The $\T^{n+m}$-equivariant homeomorphism
    \eqref{eq:zethomeomorphism} induces a graded ring isomorphism
    on equivariant cohomology which is
    $H^*(B\T^{n+m};\Z)$-linear.
    The result follows from \Cref{thm:generalpolyhedral}.
\end{proof}

\begin{remark}\label{rem:panels}
A \emph{panel structure} on a CW-complex $X$ is a finite indexed family $\{X_i\}_{i\in[m]}$ of closed subspaces, called \emph{panels}, each of which is a CW-subcomplex. Such structures were used by Davis in the study of reflection groups \cite[Section~6]{MR690848}, and Yu used them to define polyhedral products over spaces with faces \cite{yu2024panelstructures}. A nice manifold with corners is a special case: its facets $F_1, \ldots, F_m$ form a panel structure, with $X_J=F_{\cap J}$ for every $J\subseteq[m]$.

The arguments of this and the previous section use only this structure, together with the CW-decomposition of \Cref{prop:CWdecomp}. They therefore apply verbatim to any CW-complex $X$ with a panel structure, provided the polyhedral product admits a CW-structure in which the refined subproducts $\mathbf D(J, L)$ are subcomplexes. In that generality, \Cref{thm:generalpolyhedral} reads with $X_J=\bigcap_{j\in J}X_j$ in place of $F_{\cap J}$ and with the ring $\bigoplus_{J\subseteq[m]}H^*(X_J;\Z)\otimes\Z[\mathbf x_J^+]$ in place of the topological face ring $\Z\langle Q\rangle$.
\end{remark}

    \appendix
    \section{Existence of the quasi-cubical decomposition}\label{app:quasicubical}
In this appendix we prove the existence of quasi-cubical decompositions, \Cref{lem:quasicubical}, using a collar germ. We first recall the notion of a manifold with corners, then give the definition of a quasi-cubical decomposition together with the existence result, after which we prove the embedding lemma of \Cref{section:cubicalization}. The proof of the existence result is given in the last part.

\begin{definition}\label{def:nicemfdwithcorner}
    A $d$-manifold with corners $Q$ is a space locally modeled on open subsets of $\R^d_{\geq 0}$, with smooth transition maps.
For each $x=(x_1, \cdots, x_d)\in \R^d_{\geq 0}$, define $c(x) \in \mathbb{N}$ to be the number of coordinates $x_i$ that are zero. By considering local charts of $Q$, we get a well-defined function $c: Q \to \mathbb{N}$. For each $q\in Q$, we call $c(q)$ the codimension of $q$. An \emph{open face of codimension-$k$} is a connected component of $c^{-1}(k)$. A \emph{(closed) face of codimension-$k$} is the closure of an open face. A face of codimension one is called a facet of $Q$. Note that a face of codimension-$0$ is a connected component of $Q$.
For a face $f$ of $Q$, let $f^{\circ}$ denote its relative interior, which is an open face.  
\end{definition}

\subsection{Quasi-cubical decompositions and the embedding lemma}\label{subsec:app-quasicubical}
\Cref{def:quasicubical} below is the formal version of the decomposition described in \Cref{section:cubicalization}.

\begin{definition}\label{def:quasicubical}
    For $\delta >0$, a \emph{quasi-cubical decomposition} of $Q$ consists of
\begin{enumerate}
    \item A subset $Q^{\delta}$ of $Q$, which is also a nice manifold with corners.
    \item An isotopy $H_{\delta}: Q\times [0, 1] \to Q$ such that $H_{\delta}(-, 0)=\id$ and $H_{\delta}(-, 1)(Q) = Q^{\delta}$. 
    \item For $J\subseteq [m]$, let $F^{\delta}_{\cap J} := (H_{\delta}(-, 1))(F_{\cap J})$.
    \item There is an embedding $G: Q \to Q^{\delta} \times [0, \delta]^m$ whose image is 
    \begin{equation}\label{eq:imageofG}
        G(Q)=\bigcup_{J \subseteq [m]} F_{\cap J}^{\delta} \times \prod_{j\in J} [0, \delta]_{(j)} \times \prod_{j\notin J} \delta_{(j)}.
    \end{equation}
    The restriction of $G$ to $Q^{\delta}$ is the identity when $Q^{\delta}$ is identified with $Q^{\delta} \times \prod_{j\in [m]} \delta_{(j)}$. More generally, for $J\neq \varnothing$, we have 
    \begin{equation}
        G(F_{\cap J}) = \bigcup_{J \subseteq J'} F_{\cap J'}^{\delta}  \times \prod_{j\in J} 0_{(j)} \times \prod_{j \in J'\setminus J} [0, \delta]_{(j)} \times \prod_{j\notin J'} \delta_{(j)} .      
    \end{equation}
\end{enumerate}
\end{definition}

\begin{lemma}\label{lem:quasicubical}
    If $Q$ is a nice manifold with corners, it admits a quasi-cubical decomposition.  
\end{lemma}

\begin{proof}[Proof of \Cref{lemma:embeddingofQandxi}]
    Take $\delta=1$ in \Cref{lem:quasicubical} and keep the notation of the definition.
    Let $\Gamma_{\delta}: Q^{\delta} \to Q$ denote the inverse of $H_{\delta}(-, 1)$, and set $\hat f=(\Gamma_{\delta}\times\id)\circ G$.
    Since $\Gamma_{\delta}(F^\delta_{\cap J})=F_{\cap J}$, the image of this embedding and \eqref{eq:preimageoff} follow from the corresponding properties of $G$.
    We need to lift $\hat f$ to a bundle embedding $\hat f_*:\xi\to\xi\times[0, 1]^m$.

    The isotopy $(q, t)\mapsto H_{\delta}(\Gamma_{\delta}(q), t)$ provides an isotopy between $\Gamma_{\delta}$ and the inclusion $i:Q^{\delta}\to Q$.
    Let $pr_1$ denote the projection $Q \times [0, 1]^m \to Q$.
Consider the following diagram where each square is a pullback square: 
\[\begin{tikzcd}[column sep=huge]
	{( pr_1\circ (\Gamma_{\delta}\times \id )\circ G)^*(\xi)} & {\Gamma_{\delta}^* (\xi) \times [0, 1]^m} & {\xi\times [0, 1]^m} & \xi \\
	Q & {Q^{\delta} \times[0, 1]^m} & {Q\times[0, 1]^m} & Q
	\arrow["{G_*}", from=1-1, to=1-2]
	\arrow[from=1-1, to=2-1]
	\arrow["{(\Gamma_{\delta})_* \times \id}", from=1-2, to=1-3]
	\arrow[from=1-2, to=2-2]
	\arrow[from=1-3, to=1-4]
	\arrow[from=1-3, to=2-3]
	\arrow[from=1-4, to=2-4]
	\arrow["G", from=2-1, to=2-2]
	\arrow["{\Gamma_{\delta}\times \id}", from=2-2, to=2-3]
	\arrow["{pr_1}", from=2-3, to=2-4]
\end{tikzcd}.\]
Let $pr_1^{\delta}$ denote the projection $Q^{\delta} \times [0, 1]^m \to {Q}^{\delta}$. Then we have 
\(
    \Gamma_{\delta} \circ pr_1^{\delta} = pr_1\circ (\Gamma_{\delta} \times \id)
\). 
So
\begin{equation}
    ( pr_1\circ (\Gamma_{\delta}\times \id )\circ G)^*(\xi) = (\Gamma_{\delta} \circ pr_1^{\delta} \circ G)^*(\xi).
\end{equation}
Let $S(-, -): Q^{\delta} \times [0, \delta]^m \times [0, 1]\to Q^{\delta} \times [0, \delta]^m$ be a homotopy defined by
\[
    \bigl( (q, x_1, \cdots, x_m), t \bigr) \mapsto (q, t \delta + (1-t)x_1, \cdots, t\delta + (1-t)x_m). 
\]
$S$ preserves each block of $G(Q)$ and fixes $Q^\delta\times (\delta, \ldots, \delta)$. Hence we can define a homotopy on $Q$ by $G^{-1}\circ S(-, t) \circ G$. Note that
$$G^{-1}\circ S(-, 0) \circ G =\id \quad \text{and} \quad G^{-1}\circ S(-, 1) \circ G = i\circ pr_1^{\delta} \circ G. $$ Thus, we have a homotopy:
\[
    i\circ pr_1^{\delta} \circ G \simeq \id.
\]
Since $\Gamma_{\delta} \simeq i$, we have $\Gamma_{\delta}\circ pr_1^{\delta}\circ G \simeq \id$. 
By homotopy invariance of principal bundles, we have
\begin{equation*}
    ( \Gamma_{\delta}\circ pr_1^{\delta} \circ G)^*(\xi) \cong \xi. 
\end{equation*}
Choose an isomorphism from $\xi$ to this pullback bundle and compose it with the top maps of the left two pullback squares to obtain $\hat f_*$.
Since $\hat f$ is an embedding, the pullback projection identifies its source homeomorphically with $(\pi_\xi\times\id)^{-1}(\widehat Q)$.
Thus $\hat f_*$ is a $\T^m$-equivariant embedding and \eqref{cd:pullbackbundle} is a pullback square.
\end{proof}

\subsection{Proof of the existence statement}\label{subsec:app-existence}
Consider the case $m=1$. Then $Q$ is a manifold with boundary, whose only facet is $\partial Q$. The existence of the quasi-cubical decomposition follows from the collar neighborhood theorem.
In the general case, the proof depends on the \emph{collar germ} of a manifold with corners, introduced in \cite{MR226674} and defined below; it generalizes the collar neighborhood theorem to the case of several facets.

Let $\mathcal{F}=\{F_1, \cdots, F_m\}$ be the set of facets of $Q$. A collar germ $\kappa$ is a collection of collar embeddings
\begin{equation}
    \kappa_i: F_i \times [0, \infty) \to Q, \ \forall 1\leq i\leq m,
\end{equation}
such that $\kappa_i|_{F_i \times \{0\}}$ is the inclusion $F_i \subset Q$. Furthermore, for any $p \in F_i\cap F_j$ where $i\neq j$, we have 
\begin{equation}\label{eq:compatible}
    \kappa_i(\kappa_j(p, x_j), x_i) = \kappa_j(\kappa_i(p, x_i), x_j). 
\end{equation}
It follows from \cite[Section 3.4]{MR226674} that a collar germ exists for any nice manifold with corners. Intuitively, each $\kappa_i$ gives the facet $F_i$ a collar inside $Q$, and the compatibility condition says that, near the intersection of two facets, these collars commute. 

For the coordinate formulas below, relabel the facets so that $J=\{1, \ldots, s\}$.
Let $\sigma$ be a permutation on $s$ elements. Then for any $p\in F_{\cap J}$ and $(x_1, x_2, \cdots, x_s) \in \prod_{j\in J} [0, \infty)_{(j)}$, 
the compatibility condition in \eqref{eq:compatible} can be generalized as: 
\begin{equation}\label{eq:generalcompatible}
    \kappa_{1}(\kappa_{2}(\cdots \kappa_{s}(p, x_s) \cdots , x_{2}), x_{1}) 
    = \kappa_{{\sigma(1)}}(\kappa_{{\sigma(2)}}(\cdots \kappa_{{\sigma(s)}}(p, x_{\sigma(s)}) \cdots , x_{\sigma(2)}), x_{\sigma(1)}) . 
\end{equation}

The compatibility condition \cref{eq:compatible} also implies that for any $p\in F_i \cap F_j$ and $x_j \in [0, \infty)$, we have $\kappa_j(p, x_j) \in F_{i}$. 
In general, if $p\in F_1\cap F_2  \cap \cdots \cap F_s =F_{\cap J}$, then for any $j\in J$, we have $\kappa_j(p, x_j) \in F_{\cap (J\setminus\{j\})}$. 
This observation allows us to define a set of embeddings among $F_{\cap J} \times \prod_{j\in J} [0, \infty)_{(j)}$. 
Let $J'=\{1, \cdots, s, \cdots, r\}$ be a set which contains $J$ and $F_{\cap J'}\neq \varnothing$. Define
\begin{equation}\label{eq:defofembeddings}
    \begin{array}{c}
            \kappa_{(J', J)}: F_{\cap J'} \times \prod_{j\in J'} [0, \infty)_{(j)} \to F_{\cap J} \times \prod_{j\in J} [0, \infty)_{(j)} \\ [8pt]
    (p, x_1, \cdots, x_s, \cdots, x_r) \mapsto \bigl( \kappa_{s+1}( \kappa_{s+2} (\cdots \kappa_r(p, x_r) \cdots,  x_{s+2} ), x_{s+1}), x_1, \cdots, x_s \bigr). 
    \end{array}
\end{equation}
In particular, $\kappa_{(J, \varnothing)}$ is an embedding $F_{\cap J} \times \prod_{j\in J} [0, \infty) \to Q$. Moreover, the embedding $\kappa_{(\{i\}, \varnothing)}$ is exactly $\kappa_i$. 
The compatibility condition for collar germs implies for any $J\subseteq J' \subseteq J''$, we have 
\begin{equation}\label{eq:concisegeneralcompatible}
    \kappa_{(J', J)} \circ \kappa_{(J'', J')} = \kappa_{(J'', J)}. 
\end{equation}

By shrinking the image of $\kappa_i$ and reparameterizing it if necessary, we may assume for each $J\subseteq [m]$ that 
\begin{equation}\label{eq:conditionforcollargerm}
    \im \kappa_{(J, \varnothing)} = \bigcap\limits_{j \in J} \im \kappa_j. 
\end{equation}
For these collars, $\kappa_i(F_i\times[0, a])$ is closed in $Q$ for every $a>0$.
Condition \eqref{eq:conditionforcollargerm} implies that for any $J', J'' \subseteq [m]$, we have $ \im \kappa_{(J', \varnothing)}  \cap  \im \kappa_{(J'', \varnothing)} =  \im \kappa_{(J' \cup J'', \varnothing)} $. Moreover, if $J \subseteq J' \cap J''$,
combined with \eqref{eq:concisegeneralcompatible}, then we have
\begin{equation*}
    \kappa_{(J, \varnothing)}\bigl(\im \kappa_{(J', J)} \cap \im \kappa_{(J'', J)}\bigr) = \im \kappa_{(J', \varnothing)} \cap \im \kappa_{(J'', \varnothing)}  
    =\im \kappa_{(J' \cup J'', \varnothing)} = 
    \kappa_{(J, \varnothing)} \bigl(\im \kappa_{(J' \cup J'', J)}\bigr). 
\end{equation*}
Since $\kappa_{(J, \varnothing)}$ is an embedding, we can cancel it from the leftmost and rightmost sides of the above equalities and get
\begin{equation}\label{eq:overlap}
    \im \kappa_{(J', J)} \cap \im \kappa_{(J'', J)} = \im \kappa_{(J' \cup J'', J)}. 
\end{equation}
By considering the definitions of the embeddings in the above equation, for any subset $Y \subseteq \prod_{j\in J} [0, \infty)_{(j)}$, we obtain a more general intersection formula:
\begin{equation}\label{eq:generaloverlap}
    \begin{aligned}
        &\kappa_{(J', J)}\Bigl(F_{\cap J'} \times Y \times \prod_{j\in J'\setminus J} [0, \infty)_{(j)}\Bigr) \bigcap
        \kappa_{(J'', J)}\Bigl(F_{\cap J''} \times Y \times \prod_{j\in J''\setminus J} [0, \infty)_{(j)}\Bigr) \\
        &= \kappa_{(J'\cup J'', J)}\Bigl(F_{\cap (J'\cup J'')} \times Y \times \prod_{j\in (J'\cup J'')\setminus J} [0, \infty)_{(j)}\Bigr).
    \end{aligned}
\end{equation}

In \eqref{eq:defofembeddings}, if we take each $x_j=0$ for $j\in J$, then the image of $\kappa_{(J', J)}$ for such a point lies in $F_{\cap J} \times \prod_{j\in J} 0_{(j)}$. 
Hence, for any $\delta>0$, we can define a subset $\widetilde{F}_{\cap J}^{\delta}$ of $F_{\cap J}$ by 
\begin{equation*}
    \widetilde{F}_{\cap J}^{\delta} \times \prod_{j\in J} 0_{(j)} := F_{\cap J}\times \prod_{j\in J} 0_{(j)} - \bigcup_{J\subsetneq J'} \kappa_{(J', J)} \bigl( F_{\cap J'} \times \prod_{j\in J} 0_{(j)} \times \prod_{j \in J' \setminus J} [0, \delta)_{(j)}  \bigr). 
\end{equation*}
Here and in the following, we use the notation $\bigcup_{J \subsetneq J'}$ to denote the union over all $J'$ such that $J \subsetneq J' \subseteq [m]$. 
Clearly, for any $Y \subseteq \prod_{j\in J}[0, \infty)_{(j)}$, we have
\begin{equation*}
    \widetilde{F}_{\cap J}^{\delta} \times Y =
    F_{\cap J} \times Y - \bigcup_{J \subsetneq J'} \kappa_{(J', J)} \bigl(F_{\cap J'} \times Y \times \prod_{j\in J'\setminus J}[0, \delta)_{(j)} \bigr). 
\end{equation*}
The above equation can be rewritten as a decomposition of $F_{\cap J} \times Y$ into two subspaces:
\begin{equation}\label{eq:smallproduct}
    F_{\cap J} \times Y = \bigl( \widetilde{F}_{\cap J}^{\delta} \times Y \bigr) \bigcup  \Bigl( \bigcup_{J \subsetneq J'} \kappa_{(J', J)} \bigl(F_{\cap J'} \times Y \times \prod_{j\in J'\setminus J}[0, \delta)_{(j)} \bigr) \Bigr)
\end{equation}
where the two subspaces are disjoint. Also, we have
\begin{equation}\label{eq:smallproduct2}
    F_{\cap J} \times Y = \bigl( \widetilde{F}_{\cap J}^{\delta} \times Y \bigr) \bigcup  \Bigl( \bigcup_{J \subsetneq J'} \kappa_{(J', J)} \bigl(F_{\cap J'} \times Y \times \prod_{j\in J'\setminus J}[0, \delta]_{(j)} \bigr) \Bigr).
\end{equation}
Taking $J=\varnothing$, we have a decomposition of $Q$ as
\begin{equation}\label{eq:decomposeQ}
    Q = \bigcup_{J' \subseteq [m]} \kappa_{(J', \varnothing)} (\widetilde{F}_{\cap J'}^{\delta} \times \prod_{j\in J'} [0, \delta]_{(j)} ).
\end{equation}
This is a finite closed cover: $\widetilde{F}_{\cap J'}^{\delta}$ is closed in $F_{\cap J'}$, and the image of the bounded multiple collar is the intersection of the corresponding closed collar strips by \eqref{eq:conditionforcollargerm}.

\begin{proposition}\label{prop:overlap}
    For $J\subseteq J'$ and any $Y\subseteq\prod_{j\in J}[0, \infty)_{(j)}$, we have
        \begin{equation}\label{eq:intersection1}
            \begin{aligned}
                & \kappa_{(J' , J)}\bigl(F_{\cap J'} \times Y \times \prod_{j\in J' \setminus J} [0, \infty)_{(j)} \bigr) \bigcap \bigl( \widetilde{F}^{\delta}_{\cap J} \times Y \bigr) \\
                &=\kappa_{(J', J)} \bigl( \widetilde{F}^{\delta}_{\cap J'} \times Y \times \prod_{j\in J'\setminus J} [\delta, \infty)_{(j)} \bigr).
            \end{aligned}
        \end{equation}
\end{proposition}

\begin{proof}
    The case $J=J'$ is immediate. In the definition of $\widetilde{F}^{\delta}_{\cap J}$, it suffices to take the union over $J\cup\{i\}$ with $i\notin J$: by \eqref{eq:concisegeneralcompatible}, every term indexed by a larger set is contained in one of these terms.
    Thus a point is removed whenever at least one of its collar coordinates normal to a facet of $F_{\cap J}$ is less than $\delta$.

    Consider $\kappa_{(J', J)}(p, x)$ with $p\in F_{\cap J'}$ and $x|_J\in Y$.
    For $i\in J'\setminus J$, the corresponding collar coordinate is $x_i$.
    For $i\notin J'$, \eqref{eq:overlap} and \eqref{eq:concisegeneralcompatible} show that this point lies in the collar of $F_{\cap(J\cup\{i\})}$ precisely when $p$ lies in the collar of $F_{\cap(J'\cup\{i\})}$ in $F_{\cap J'}$, with the same $i$-coordinate.
    Hence $\kappa_{(J', J)}(p, x)$ lies in $\widetilde{F}^{\delta}_{\cap J}\times Y$ if and only if $p\in\widetilde{F}^{\delta}_{\cap J'}$ and $x_i\geq\delta$ for every $i\in J'\setminus J$. This proves \eqref{eq:intersection1}.
\end{proof}

\begin{corollary}\label{coro:intersections}
    Let $J\subseteq J'\cap J''$, and put $U=J'\cup J''$ and $V=J'\cap J''$.
    For any $Y\subseteq\prod_{j\in J}[0, \infty)_{(j)}$, we have
    \begin{equation}
        \begin{aligned}
            &\kappa_{(J', J)}\Bigl(\widetilde{F}^{\delta}_{\cap J'}\times Y\times\prod_{j\in J'\setminus J}[0, \infty)_{(j)}\Bigr) \\
            &\qquad\cap\kappa_{(J'', J)}\Bigl(\widetilde{F}^{\delta}_{\cap J''}\times Y\times\prod_{j\in J''\setminus J}[0, \infty)_{(j)}\Bigr) \\
            &=\kappa_{(U, J)}\Bigl(\widetilde{F}^{\delta}_{\cap U}\times Y
            \times\prod_{j\in V\setminus J}[0, \infty)_{(j)}
            \times\prod_{j\in U\setminus V}[\delta, \infty)_{(j)}\Bigr).
        \end{aligned}
    \end{equation}
\end{corollary}

\begin{proof}
    By \eqref{eq:generaloverlap}, every point in the intersection has a unique expression $\kappa_{(U, J)}(p, x)$ with $p\in F_{\cap U}$ and $x|_J\in Y$.
    By \Cref{prop:overlap}, it belongs to the first set precisely when $p\in\widetilde{F}^{\delta}_{\cap U}$ and $x_j\geq\delta$ for $j\in U\setminus J'$.
    Membership in the second set similarly requires $x_j\geq\delta$ for $j\in U\setminus J''$.
    Together these conditions impose $x_j\geq\delta$ exactly for $j\in U\setminus V$, leaving the coordinates in $V\setminus J$ unrestricted.
\end{proof}

\begin{corollary}\label{coro:cubeintersections1}
    Let $J\subseteq J'\cap J''$, and put $U=J'\cup J''$ and $V=J'\cap J''$. Then
    \begin{equation}
        \begin{aligned}
            &\kappa_{(J', J)}\Bigl(\widetilde{F}^{\delta}_{\cap J'}\times\prod_{j\in J'}[0, \delta]_{(j)}\Bigr)
            \cap\kappa_{(J'', J)}\Bigl(\widetilde{F}^{\delta}_{\cap J''}\times\prod_{j\in J''}[0, \delta]_{(j)}\Bigr) \\
            &=\kappa_{(U, J)}\Bigl(\widetilde{F}^{\delta}_{\cap U}
            \times\prod_{j\in V}[0, \delta]_{(j)}
            \times\prod_{j\in U\setminus V}\{\delta\}_{(j)}\Bigr),
        \end{aligned}
    \end{equation}
    and
    \begin{equation}\label{eq:consistentofF}
        \begin{aligned}
            &\kappa_{(J', J)}\Bigl(\widetilde{F}^{\delta}_{\cap J'}\times\prod_{j\in J'}\{\delta\}_{(j)}\Bigr)
            \cap\kappa_{(J'', J)}\Bigl(\widetilde{F}^{\delta}_{\cap J''}\times\prod_{j\in J''}\{\delta\}_{(j)}\Bigr) \\
            &=\kappa_{(U, J)}\Bigl(\widetilde{F}^{\delta}_{\cap U}\times\prod_{j\in U}\{\delta\}_{(j)}\Bigr).
        \end{aligned}
    \end{equation}
\end{corollary}

\begin{proof}
    In \Cref{coro:intersections}, restricting all collar coordinates to $[0, \delta]$ forces those in $U\setminus V$ to equal $\delta$, while those in $V$ remain free in $[0, \delta]$. Restricting the remaining coordinates to $\delta$ gives \eqref{eq:consistentofF}.
\end{proof}
In particular, taking $J=\varnothing$ gives the intersections of the blocks in \eqref{eq:decomposeQ}.

\begin{proposition}
    For a given $\delta >0$, there exists an isotopy
\[
    H_{(J, \delta)}(-, -): F_{\cap J} \times \prod_{j\in J}[0, \infty)_{(j)} \times [0, 1] \to F_{\cap J} \times \prod_{j\in J}[0, \infty)_{(j)}
\]
such that 
\begin{enumerate}
    \item $H_{(J, \delta)}( -, 0)= \id$, $H_{(J, \delta)}( -, 1)$ maps $F_{\cap J} \times \prod_{j\in J} [0, \infty)_{(j)}$ onto $\widetilde{F}^{\delta }_{\cap J} \times \prod_{j\in J} [\delta, \infty)_{(j)}$ and
    maps $F_{\cap J} \times \prod_{j\in J} 0_{(j)}$ onto $\widetilde{F}^{\delta }_{\cap J} \times \prod_{j\in J} \delta_{(j)}$. 
    \item for any $t\in [0, 1]$, the restriction $H_{(J, \delta)}(- , t)$ on  $\widetilde{F}^{2\delta}_{\cap J} \times \bigl( \prod_{j\in J} [0, \infty)_{j}\bigr)$ is given by
    \begin{equation}\label{eq:defofisotopyonR}
        (p, x_1, x_2, \cdots, x_s) \mapsto (p, h_{\delta}(x_1, t), \cdots, h_{\delta}(x_s, t))
    \end{equation}
    where $J=\{1, \ldots, s\}$ after relabelling and $h_\delta$ is the isotopy on $[0, \infty)$ defined by
    \[
        h_{\delta} (x, t) = 
        \begin{cases}
            (1-t)x+t\dfrac{x+2\delta}{2}, &x\leq2\delta, \\
            x, &x\geq2\delta.
        \end{cases}
    \]
    \item These isotopies are compatible in the sense that for any $J\subseteq J'$, we have
    \begin{equation}
        \kappa_{(J', J)}\circ H_{(J', \delta)} = H_{(J, \delta)} \circ \kappa_{(J', J)}. 
    \end{equation}
\end{enumerate}
\end{proposition}

\begin{proof}
    We construct $H_{(J, \delta)}$ by induction on the cardinality of $J$ from $m$ down to $0$. 

    First, assume that $J$ is maximal among the subsets of $[m]$ such that $F_{\cap J} \neq \varnothing$. In this case, $\widetilde{F}^{\delta}_{\cap J} = F_{\cap J}$. 
    Define $H_{(J, \delta)}$ by \eqref{eq:defofisotopyonR}. The first two assertions follow directly, and the third is vacuous by maximality.
    
    Suppose we have constructed $H_{(J', \delta)}$ for any $|J'|>s$. We construct $H_{(J, \delta)}$ for $|J|=s$. Define $H_{(J, \delta)}$ by 
    \begin{equation*}
        H_{(J, \delta)} ((p, x), t) = 
        \begin{cases}
            (p, h_{\delta}(x_1, t), \cdots, h_{\delta}(x_s, t)), \ \ p \in \widetilde{F}^{2 \delta}_{\cap J} ; \\
            \kappa_{(J', J)} \circ H_{(J', \delta)} \circ \kappa_{(J', J)}^{-1}(p, x), \ \ (p, x) \in \im \kappa_{(J', J)}. 
        \end{cases}
    \end{equation*}
    The well-definedness of $H_{(J, \delta)}$ needs to be verified:
    \begin{itemize}
        \item If $(p, x) \in \im \kappa_{(J', J)} \cap \im \kappa_{(J'', J)}$, since $\kappa_{(J', J)}\circ \kappa_{(J'\cup J'', J')} = \kappa_{(J'', J)}\circ \kappa_{(J'\cup J'', J'')}$, we have
        \begin{equation*}
            \begin{aligned}
                &\kappa_{(J', J)} \circ H_{(J', \delta)} \circ \kappa_{(J', J)}^{-1}(p, x)  \\
            = &\kappa_{(J', J)} \circ \kappa_{(J'\cup J'', J')} \circ H_{(J'\cup J'', \delta)} \circ \kappa_{(J'\cup J'', J')}^{-1} \circ \kappa_{(J', J)}^{-1}(p, x) \\ 
            = &\kappa_{(J'', J)} \circ \kappa_{(J'\cup J'', J'')} \circ H_{(J'\cup J'', \delta)} \circ \kappa_{(J'\cup J'', J'')}^{-1} \circ \kappa_{(J'', J)}^{-1}(p, x) \\
            = &\kappa_{(J'', J)} \circ H_{(J'', \delta)} \circ \kappa_{(J'', J)}^{-1}(p, x).
            \end{aligned}
        \end{equation*}
        The first equality uses the following equation by induction:
        \[ 
            H_{(J', \delta)} \circ \kappa_{(J'\cup J'', J')} =\kappa_{(J'\cup J'', J')} \circ H_{(J'\cup J'', \delta)}
        \]
        
        \item If $(p, x) \in \im \kappa_{(J', J)} \cap \widetilde{F}^{2 \delta}_{\cap J} \times \prod_{j\in J} [0, \infty)_{(j)}$, it follows from \Cref{prop:overlap} that
        \[ 
            \im \kappa_{(J', J)} \bigcap \widetilde{F}^{2 \delta}_{\cap J} \times \prod_{j\in J} [0, \infty)_{(j)} = \kappa_{(J', J)} \bigl( \widetilde{F}^{2 \delta}_{\cap J'} \times \prod_{j\in J} [0, \infty)_{(j)} \times \prod_{j\in J'\setminus J} [2\delta, \infty) \bigr). 
        \]
        Relabel so that $J=\{1, \ldots, s\}$ and $J'=\{1, \ldots, r\}$. Let $(q, x_1, \cdots, x_r) \in \widetilde{F}^{2 \delta}_{\cap J'} \times \prod_{j\in J} [0, \infty)_{(j)} \times \prod_{j\in J'\setminus J} [2\delta, \infty)_{(j)}$ such that
        \[p = \kappa_{s+1}( \kappa_{s+2}( \cdots \kappa_r(q, x_r) \cdots , x_{s+2}), x_{s+1} ) \in \widetilde{F}^{\delta}_{\cap J}. \] 
        Note that $h_{\delta}( - , t) = \id $ on $[2\delta, \infty)$ for any $t\in [0, 1]$, hence we have 
        \[
        \begin{aligned}
            &\kappa_{(J', J)} \circ H_{(J', \delta)} (q, x_1, x_2, \cdots, x_r) \\ 
            = & \kappa_{(J', J)} (q, h_{\delta}(x_1, t), \cdots, h_{\delta}(x_s, t), x_{s+1}, \cdots, x_r) \\
            = & \bigl( \kappa_{s+1}( \kappa_{s+2}( \cdots \kappa_r(q, x_r) \cdots , x_{s+2}), x_{s+1}), h_{\delta}(x_1, t), \cdots, h_{\delta}(x_s, t) \bigr) . 
        \end{aligned}\]
        Hence, we have 
        \[
        \kappa_{(J', J)} \circ H_{(J', \delta)} \circ \kappa_{(J', J)}^{-1}(p, x) = (p, h_{\delta}(x_1, t), \cdots, h_{\delta}(x_s, t)). 
        \]
    \end{itemize}
    To verify continuity and the embedding property, apply the decomposition in \eqref{eq:decomposeQ} to $F_{\cap J}$ with $2\delta$ in place of $\delta$, and take its product with $\prod_{j\in J}[0, \infty)_{(j)}$.
    This gives a finite closed cover. In the collar coordinates of each block, $H_{(J, \delta)}$ fixes the base point and applies $h_\delta$ to each interval coordinate, so it is jointly continuous on the block and $[0, 1]$.
    For each $t$, the block map is a homeomorphism onto a closed image block, with inverse given in each interval coordinate by
    \[
        h_\delta^{-1}(r, t)=
        \begin{cases}
            \dfrac{r-t\delta}{1-t/2}, &t\delta\leq r\leq2\delta, \\
            r, &r\geq2\delta.
        \end{cases}
    \]
    By \Cref{coro:intersections}, the coordinates fixed to $2\delta$ on block intersections remain fixed under both maps. Thus the block maps and their inverses agree on overlaps.
    The pasting lemma for finite closed covers shows that $H_{(J, \delta)}$ is continuous and each $H_{(J, \delta)}(-, t)$ is an embedding.

    Properties (2) and (3) follow from the construction. It remains to verify (1) at $t=1$.
    The image of $F_{\cap J} \times \prod_{j\in J} [0, \infty)_{(j)}$ under $H_{(J, \delta)}(-, 1)$ is the union of the following two spaces:
    \begin{equation*}
        \mathbf{S}_1:= \widetilde{F}^{2\delta}_{\cap J} \times \prod_{j\in J} [\delta, \infty)_{(j)}, \quad \quad \mathbf{S}_2= \bigcup_{J\subsetneq J'} \kappa_{(J', J)} ( \widetilde{F}^{\delta}_{\cap J'} \times \prod_{j\in J'} [\delta, \infty)_{(j)} ). 
    \end{equation*}
    We need to verify that $\mathbf{S}_1 \cup \mathbf{S}_2$ is exactly
    \begin{equation*}
        \mathbf{S}:= \widetilde{F}^{\delta}_{\cap J} \times \prod_{j\in J} [\delta, \infty)_{(j)}. 
    \end{equation*}
    On the one hand, since $\mathbf{S}_1$ and $\mathbf{S}_2$ are both contained in $\mathbf{S}$, we have $\mathbf{S}_1 \cup \mathbf{S}_2 \subseteq \mathbf{S}$. 
    On the other hand, the following space is contained in $\mathbf{S}_1$
    \begin{equation*}
        F_{\cap J} \times \prod_{j\in J} [\delta, \infty)_{(j)} - \overbrace{\bigcup_{J\subsetneq J'} \kappa_{(J', J)} \bigl( F_{\cap J'} \times \prod_{j\in J} [\delta, \infty)_{(j)} \times \prod_{j\in J'\setminus J} [0, \infty)_{(j)} \bigr)}^{\mathbf{S_3}}. 
    \end{equation*}
    In \cref{eq:intersection1}, take $Y=\prod_{j\in J} [\delta, \infty)_{(j)}$
    and take the union on both sides with respect to all $J'$ which strictly contains $J$, we have 
\begin{equation}\label{eq:finialoverlap}
    \begin{aligned}
        &\Bigl( \bigcup_{J\subsetneq J'}\kappa_{(J' , J)}\bigl(F_{\cap J'} \times \prod_{j\in J} [\delta, \infty)_{(j)} \times  \prod_{j\in J' \setminus J} [0, \infty)_{(j)} \bigr) \Bigr) \bigcap \bigl( \widetilde{F}^{\delta}_{\cap J} \times \prod_{j\in J} [\delta, \infty)_{(j)}  \bigr) \\
    =  & \bigcup_{J\subsetneq J'} \kappa_{(J', J)} \bigl( \widetilde{F}^{\delta}_{\cap J'} \times \prod_{j\in J'} [\delta, \infty)_{(j)}  \bigr). 
    \end{aligned}
\end{equation}
    That is $\mathbf{S_3} \cap \mathbf{S} = \mathbf{S}_2$. Hence, we have 
    \begin{equation*}
        \begin{aligned}
            \mathbf{S}_1 \cup \mathbf{S}_2 &\supseteq (F_{\cap J} \times \prod_{j\in J} [\delta, \infty)_{(j)} - \mathbf{S}_3 ) \cup \mathbf{S}_2 
            \supseteq (\mathbf{S} -\mathbf{S}_3) \cup \mathbf{S}_2 \\
             &= (\mathbf{S} -(\mathbf{S}_3\cap \mathbf{S})) \cup \mathbf{S}_2 
             = \mathbf{S}. 
        \end{aligned}
    \end{equation*}
    Hence, we have proved that $H_{(J, \delta)}( -, 1)$ maps $F_{\cap J} \times \prod_{j\in J} [0, \infty)_{(j)}$ onto $\widetilde{F}^{\delta }_{\cap J} \times \prod_{j\in J} [\delta, \infty)_{(j)}$. 

    Since $h_\delta(x, 1)=\delta$ if and only if $x=0$, the established surjectivity also shows that $H_{(J, \delta)}(-, 1)$ maps the zero-coordinate subspace onto $\widetilde F^\delta_{\cap J}\times\prod_{j\in J}\{\delta\}_{(j)}$. This completes the induction.
\end{proof}

Set $Q^\delta=\widetilde{F}^{\delta}_{\cap\varnothing}$. The translated collar coordinates $x_i-\delta$ give $Q^\delta$ the structure of a nice manifold with corners.
Now we get a homeomorphism $H_{(\varnothing, \delta)}( - , 1) :Q \to Q^\delta$. For each $J\subseteq [m]$, we have
\begin{equation*}
    \begin{aligned}
        H_{(\varnothing, \delta)}(F_{\cap J}, 1) = & H_{(\varnothing, \delta)}(-, 1) \circ \kappa_{(J, \varnothing)}(F_{\cap J} \times \prod_{j\in J} 0_{(j)}) \\ 
    =& \kappa_{(J, \varnothing)} \circ H_{(J, \delta)}(F_{\cap J} \times \prod_{j\in J} 0_{(j)}, 1) \\ 
    =& \kappa_{(J, \varnothing)} (\widetilde{F}^{\delta}_{\cap J} \times \prod_{j\in J} \delta_{(j)}).
    \end{aligned}
\end{equation*}
Let $F^{\delta}_{\cap J}= \kappa_{(J, \varnothing)} (\widetilde{F}^{\delta}_{\cap J} \times \prod_{j\in J} \delta_{(j)})$ and $\Gamma^{\delta}$ be the inverse of 
$H_{(\varnothing, \delta)}( - , 1)$. 
Naturally, there is a series of homeomorphisms
\begin{equation*}
    (*): \widetilde{F}^{\delta}_{\cap J} \cong \widetilde{F}^{\delta}_{\cap J} \times \prod_{j\in J} \delta_{(j)} \cong F^{\delta}_{\cap J}. 
\end{equation*}
By the above homeomorphism between $F^{\delta}_{\cap J}$ and $\widetilde{F}^{\delta}_{\cap J}$, we define
\begin{equation*}
    (**): \widetilde{F}^{\delta}_{\cap J} \times \prod_{j\in J}[0, \delta]_{(j)}  {\longrightarrow} F^{\delta}_{\cap J} \times \prod_{j\in J}[0, \delta]_{(j)} \times \prod_{j\notin J} \delta_{(j)} 
\end{equation*}
where $(**)$ is defined as the identity map on the term $\prod_{j\in J}[0, \delta]_{(j)}$ and as the composition of homeomorphism on $F^{\delta}_{\cap J}$ given by 
$
      \widetilde{F}^{\delta}_{\cap J} \overset{*}{\cong}   F^{\delta}_{\cap J} \cong F^{\delta}_{\cap J} \times \prod_{j\notin J} \delta_{(j)} 
$.
According to \eqref{eq:decomposeQ}, we have 
$$Q=\bigcup_{J\subseteq [m]} \kappa_{(J, \varnothing)}( \widetilde{F}^{\delta}_{\cap J} \times \prod_{j\in J}{[0, \delta]_{(j)}}).$$
Define $G: Q \to Q^{\delta} \times [0, \delta]^m$ by defining its restriction to $\kappa_{(J, \varnothing)}( \widetilde{F}^{\delta}_{\cap J} \times \prod_{j\in J}{[0, \delta]_{(j)}})$ as the composition
\begin{equation}
    \kappa_{(J, \varnothing)}( \widetilde{F}^{\delta}_{\cap J} \times \prod_{j\in J}{[0, \delta]_{(j)}})  \longrightarrow  \widetilde{F}^{\delta}_{\cap J} \times \prod_{j\in J}{[0, \delta]_{(j)}} \longrightarrow
    {F}^{\delta}_{\cap J} \times \prod_{j\in J}{[0, \delta]_{(j)}} \times \prod_{j\notin J}{\delta_{(j)}}
\end{equation}
To check that these maps agree on overlaps, put $U=J'\cup J''$ and $V=J'\cap J''$.
By \Cref{coro:cubeintersections1}, a point in the overlap of the $J'$- and $J''$-blocks has the form $\kappa_{(U, \varnothing)}(p, x)$, where $p\in\widetilde{F}^{\delta}_{\cap U}$, $x_j\in[0, \delta]$ for $j\in V$, and $x_j=\delta$ for $j\in U\setminus V$.
By \eqref{eq:concisegeneralcompatible}, both definitions send this point to
\[
    \bigl(\kappa_{(U, \varnothing)}(p, (\delta)_{j\in U}), z\bigr),
    \qquad
    z_j=\begin{cases}x_j, &j\in U, \\ \delta, &j\notin U.\end{cases}
\]
Thus $G$ is well-defined.
Conversely, by \eqref{eq:consistentofF}, a point shared by the two image blocks has the same form, and both block inverses send it to $\kappa_{(U, \varnothing)}(p, (z_j)_{j\in U})$.
Hence $G$ maps $Q$ bijectively onto the set displayed in \eqref{eq:imageofG}.
The source blocks form a finite closed cover by \eqref{eq:decomposeQ}. The image blocks are also closed, since each $F^\delta_{\cap J}$ is the image of the closed set $F_{\cap J}$ under the homeomorphism $H_{(\varnothing, \delta)}(-, 1):Q\to Q^\delta$.
The block maps and their inverses therefore paste to continuous maps, so $G$ is an embedding.
Its restriction to $Q^\delta$ is the identity under the stated identification, and its $i$-th interval coordinate is zero precisely on $F_i$. Thus it has all the required properties.

    \printbibliography
\end{document}